\documentstyle[11pt,amssymb,amscd,epsf]{article}
\newtheorem{thm}{Theorem}[section]
\newtheorem{lem}[thm]{Lemma}
\newtheorem{prop}[thm]{Proposition}

\newtheorem{cor}[thm]{Corollary}

\newcommand{\inv}{^{-1}}
\newcommand{\iso}{\stackrel{\sim}{\longrightarrow}}
\newcommand{\tim}{^\times}

\newcommand{\Z}{{\mathbb{Z}}}
\newcommand{\R}{{\mathbb{R}}}
\newcommand{\Q}{{\mathbb{Q}}}

\newcommand{\Gm}{{\mathbb{G}}_m}

\newcommand{\cL}{{\cal{L}}}

\newcommand{\ov}{\overline}

\newcommand{\Ad}{\mbox{Ad}}

\newcommand{\om}{\Omega}

\newcommand{\bg}{{\bf  G}}
\newcommand{\bn}{{\bf N}}
\newcommand{\bt}{{\bf T}}
\newcommand{\bu}{{\bf U}}

\newcommand{\n}{\underline{n}}
\newcommand{\La}{\Lambda}
\newcommand{\Labar}{\overline{\Lambda}}
\newcommand{\Uaom}{U_{a,\Omega}}
\newcommand{\Uaomp}{U_{a,\Omega'}}

\newcommand{\Rinfty}{\R_{\geq 0, \infty}}
\begin{document}

\title{Compactification of the Bruhat-Tits building of PGL by lattices of smaller rank}
{ \author{Annette Werner \\ \small Mathematisches Institut, Universit\"at M\"unster, Einsteinstr. 62, D -  48149 M\"unster\\ \small e-mail: werner@math.uni-muenster.de}}
\date{}

\maketitle
\centerline{\bf Abstract} 
\small

In this paper we construct a compactification of the Bruhat-Tits building associated to the group $PGL(V)$ by attaching all the Bruhat-Tits buildings of $PGL(W)$ for non-trivial subspaces $W$ of $V$ as a boundary. 

\normalsize
\vspace*{0.5cm}
\centerline{{\bf MSC} (2000): 20E42, 20G25}

\section{Introduction}
Let $K$ be a finite extension of ${\mathbb{Q}}_p$ and $V$  a finite-dimensional vector space over  $K$. In this paper we construct a compactification of the Bruhat-Tits building $X$ associated to the group $PGL(V)$ by attaching all the Bruhat-Tits buildings of  $PGL(W)$ for non-trivial subspaces $W$ of $V$ as a boundary. 
 Since the vertices of these buildings correspond to homothety classes of lattices in $W$, we can also view this process as attaching to $X$ (whose underlying  simplicial complex is defined by lattices of full rank in $V$) all the lattices in $V$ of smaller rank.  

This compactification differs from both the Borel-Serre compactification and Landvogt's polyhedral compactification of $X$.  The different features of these three constructions can be illustrated in the case of a three-dimensional $V$ by looking at the compactification of one appartment: 

\epsfxsize=11cm
\epsfysize=3cm
\epsfbox{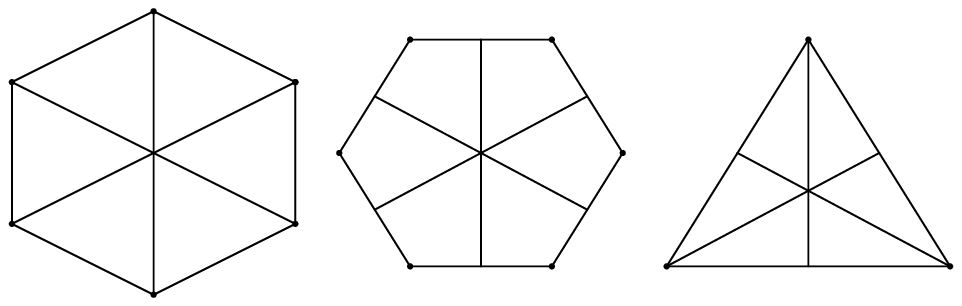}
\nopagebreak[4]
\hspace{3ex} Borel-Serre  \hspace{9ex} Landvogt \hspace{6ex} Our compactification

\pagebreak
In the Borel-Serre compactification, the points on the boundary correspond to the rays, i.e. to the half-lines emanating from the origin, and parallel lines have the same limits. In Landvogt's compactification, parallel lines have different limit points, whereas the rays in one segment (chamber) all converge to the corner vertex. In our compactification all rays contained in the two chambers around a boundary vertex converge to this vertex, and so do lines which are parallel to the middle axis. The two rays at the boundary of this double chamber (they look shorter in our picture) converge to points on the boundary lines, and their parallels converge to different points on these lines.

The idea to attach lattices of smaller rank to compactify $X$ already appeared in Mustafin's paper \cite{mu}. The goal of this paper is a generalization of Mumford's $p$-adic Schottky uniformization to higher dimensions. 
Mustafin's construction and investigation of the compactified building take up about one page. He works with lattices and defines the compactification as the union of $X$ and all the lattices in $V$ of smaller rank, i.e. he only uses a set of vertices as the boundary. His construction remains rather obscure 
(at least to the author), and does not include proofs. 

The construction of our compactification is based on the same idea of attaching lattices of smaller rank, but is entirely different. First we compactify one appartment $\La$ in $X$ (corresponding to a maximal torus $T$ in $PGL(V)$) by attaching some appartments of lower dimension corresponding to certain tori which are quotients of $T$. We define a continuous action of  the normalizer $N$ of $T$ on this compactification $\Labar$. Then we glue all the compactified appartments corresponding to maximal tori in $PGL(V)$ together. To be precise,  we take a certain compact subgroup $U_0^\wedge$ in $G$, and we define for each $x \in \Labar$ a subgroup $P_x$ of $G$, which turns out later to be the stabilizer of $x$. Our compactification $\ov{X}$ is defined as  the quotient of $U_0^\wedge \times \Labar$ by the following  equivalence relation: $(g,x) \sim(h,y)$ iff there exists some $n$ in $N$ such that $nx=y$ and $g\inv h n$ lies in $P_x$. This is similar to the construction of the building $X$. 

Then we prove the following results: $X$ is an open, dense subset of $\ov{X}$, and $\ov{X}$ carries a $G$-action compatible with the one on $X$. Besides, $\ov{X}$ is compact and contractible and can be identified with the union of all Bruhat-Tits buildings corresponding to non-zero subspaces $W$ of $V$.

In order to prove these results, we have to investigate in detail the structure of our stabilizer groups $P_x$. In particular, we prove a mixed Bruhat decomposition theorem for them. 

Of course, it would be desirable to construct similar compactifications also for buildings corresponding to other reductive groups. Besides, it would be interesting to see if there exists an analogue of our construction in the archimedean world of symmetric spaces.

{\bf Acknowledgements: } I would like to thank C. Deninger, E. Landvogt and P. Schneider for helpful discussions concerning this paper. 

\section{The Bruhat-Tits building for $PGL$}

Throughout this paper we denote by $K$ a finite extension of $\Q_p$, by $R$ its valuation ring and by $k$ the residue class field. Besides, $v$ is the valuation map, normalized so that it maps a prime element to $1$.

We adopt the convention that ``$\subset$'' always means strict subset, whereas we write ``$\subseteq$'', if equality is permitted.

Let $V$ be an $n$-dimensional vector space over $K$. Let us recall the definition of the Bruhat-Tits building for the group $\bg = PGL(V)$ (see \cite{brti1} and \cite{la}).

We fix a maximal $K$-split torus $\bt$ and let $\bn = N_G \bt$ be its normalizer. Note that $\bt$ is equal to its centralizer in $\bg$. We write  $G = \bg(K)$, $T=\bt(K)$ and $N= \bn(K)$ for the groups of rational points. 
By $X_\ast(\bt)$ respectively $X^\ast(\bt)$ we denote the cocharacter respectively the character group of $\bt$. 
We have a natural perfect pairing
\begin{eqnarray*}
<,>: & X_\ast(\bt) \times X^\ast(\bt) & \longrightarrow  \Z \\
~ & (\lambda, \chi) & \longmapsto  <\lambda, \chi> ,
\end{eqnarray*}
where $<\lambda, \chi>$ is the integer such that $\chi\circ \lambda(t) = t^{<\lambda, \chi>}$ for all $t \in \Gm$.
Let $\Lambda$ be the $\R$-vector space $\Lambda= X_\ast(\bt) \otimes_\Z \R$. We can identify the dual space $\Lambda^\ast$ with $X^\ast(\bt) \otimes_\Z \R$, and extend $<,>$ to a pairing
\[ <,>: \Lambda \times \Lambda^\ast \longrightarrow \R.\]
Since $<,>$ is perfect, there exists a unique homomorphism $\nu: T \rightarrow \Lambda$ such that
\[ < \nu(z), \chi> = - v(\chi(z))\]
for all $z \in T$ and $\chi \in X^\ast(\bt)$ (compare \cite{la}, Lemma 1.1). 
Besides, by \cite{la}, Proposition 1.8, there exists an  affine $\Lambda$-space $A$ together with a homomorphism $\nu: N \rightarrow \mbox{Aff}(A)$ extending $\nu: T \rightarrow \Lambda$. Here $\mbox{Aff}(A)$ denotes the space of affine bijections $A \rightarrow A$. The pair $(A, \nu)$ is unique up to unique isomorphism. It is called the empty appartment defined by $\bt$.

Let $g$ be the Lie algebra of $\bg$. We have the root decomposition
\[ g = g^T \oplus \bigoplus_{a \in \Phi} g_a,\]
where  $\Phi = \Phi(\bt,\bg)$ is the set of roots and where $g^T = \{X \in g: \Ad(t) X = X \mbox{ for all } t \in T\}$ and $g_a = \{ X \in g: \Ad(t)X = a(t) X \mbox{ for all }t \in T\}$ for all $a  \in \Phi$ (see \cite{bo}, 8.17 and 21.1). By \cite{bo}, 21.6, $\Phi$ is a root system in $\Lambda^\ast$ with Weyl group $W = N/T$. 
For all $a\in \Phi$ there exists a unique closed, connected, unipotent subgroup $\bu_a$ of $\bg$ which is normalized by $\bt$ and has Lie algebra $g_a$ (see \cite{bo}, 21.9). We denote the $K$-rational points of $\bu_a$ by $U_a$. 

In our case $\bg = PGL(V)$ we can describe these data explicitely. Our torus $\bt$ is the image of a maximal split torus $\bt^\sim$ in $GL(V)$. Hence there exists a basis $v_1, \ldots ,v_n$ of $V$ such that 
$ \bt^\sim$ is the group of diagonal matrices in $GL(V)$ with respect to  $v_1,\ldots,v_n$. From now on we will fix such a basis. 
Let $\bn^\sim$ be the normalizer of $\bt^\sim$ in $GL(V)$. Then $\bn$ is the image of $\bn^\sim$ in $PGL(V)$ by \cite {bo}, 22.6. Hence $N$ is the semidirect product of  $ T$ and the group of permutation matrices, which is isomorphic to $W = N/T$.  

Since $W$ is the Weyl group corresponding to $\Phi$, it acts as a group of reflections on $\Lambda$, and we have a natural homomorphism
\[ W \longrightarrow GL(\Lambda).\]
Since $\mbox{Aff}(\Lambda) = \Lambda \rtimes GL(\Lambda)$, we can use this map together with $\nu: T \rightarrow \Lambda$ to define
\[ \nu: N = T \rtimes W \longrightarrow \Lambda \rtimes GL(\Lambda) = \mbox{Aff} (\Lambda).\]
Hence $(\Lambda, \nu)$ is an empty appartment, and we write from now on $A = \Lambda$.

Denote by $\chi_i$ the character
\begin{eqnarray*}
\chi_i : & \bt^\sim & \longrightarrow \Gm \\
~ & \left( \begin{array}{ccc} t_1& ~ & ~ \\
~ & \ddots & ~\\
~ & ~ & t_n \end{array} \right)
& \longmapsto t_i.
\end{eqnarray*}
Then  for all $i$ and $j$  we have characters $a_{ij}:= \chi_i - \chi_j$, and 
\[ \Phi = \{ a_{ij}: i \neq j\} .\]
For $a = a_{ij}$ we define now $\bu_a^\sim$ as the subgroup of $GL(V)$ such that $\bu_a^\sim (\ov{K})$ is the group of matrices $U = (u_{kl})_{k,l}$ such that the diagonal elements $u_{kk}$ are equal to one, $u_{ij}$ is an element in $\ov{K}$ and the rest of the entries $u_{kl}$ is zero. 
Its image in $PGL(V)$ is isomorphic to $\bu^\sim_a$ and coincides with the group $\bu_a$ (see \cite{bo}, 22.6). Define
\[\psi_a : U_a \longrightarrow \Z \cup \{\infty\}\]
by mapping the matrix $U=(u_{kl})_{k,l}$ to $v(u_{ij})$. Then we put for all $l \in \Z$
\[U_{a,l} = \{ u \in U_a: \psi_a(u) \geq l \}. \]
We also define $U_{a,\infty} = \{1\}$, and $U_{a, -\infty} = U_a$. 
An affine function $\theta: \La \rightarrow \R$ of the form $\theta(x) = a(x) +l$ for some $a \in \Phi$ and some $l \in \Z$ is called an affine root. We can define an equivalence relation $\sim$ on $\La$ as follows:
\begin{eqnarray*}
x \sim y & \mbox{iff }& \theta(x) \mbox{ and } \theta(y) \mbox{ have the same sign } \\
~ & ~ & \mbox{or are both equal to } 0 \mbox{ for all affine roots } \theta.
\end{eqnarray*}
The equivalence classes with respect to this relation are called the faces of $\La$. These faces are simplices which partition $\La$ (see \cite{bou}, V, 3.9).  There exists a $W$-invariant scalar product on $\Lambda$ (uniquely determined up to scalar factor), see \cite{bou}, VI, 1.1 and 1.2, and all the reflections at affine hyperplanes are contained in $\nu(N)$ by \cite {la}, Proposition 11.8.

For all $x \in \Lambda$ let $U_x$ be the group generated by $U_{a, -a(x)}= \{u \in U_a: \psi_a(u) \geq -a(x)\}$ for all $a \in \Phi$. Besides, put $N_x= \{ n \in N: \nu(n)x=x \}$, and 
\[P_x = U_xN_x = N_x U_x.\]
Now we are ready to define the building $X= X(PGL(V))$ as
\[ X = G \times \La /\sim,\]
where the equivalence relation $\sim$ is defined as follows (see \cite{la}, 13.1):
\begin{eqnarray*}
(g,x) \sim (h,y) & \mbox{iff  there exists an element } n \in N \\
~ & \mbox{such that } \nu(n)x= y \mbox{ and } g\inv h n \in P_x.\end{eqnarray*}
We have a natural action of $G$ on $X$ via left multiplication on the first factor. The $G$-action on $X$ continues the $N$-action on $\La$, so that we will write $nx$ for our old $\nu(n) x$ if $x \in \La$ and $n \in N$.  
Besides, we can embed the appartment $\La$ in $X$, mapping $a \in \La$ to the class of $(1,a)$. 
(This is injective, see \cite{la}, Lemma 13.2.) 
For $x \in \La$ the group $P_x$ is the stabilizer of $x$. A subset of $X$ of the form $g \La$ for some $g \in G$ is called appartment in $X$. Similarly, we define the faces in $g \La$ as the subsets $g F$, where $F$ is a face in $\La$. Then two points (and even two faces) in $X$ are always contained in a common appartment (\cite{la}, Proposition 13.12 and \cite{brti1}, 7.4.18). Any appartment which contains a point of a face contains the whole face, and even its closure (see 
\cite{la}, 13.10, 13.11, and \cite{brti1}, 7.4.13, 7.4.14). We fix once and for all a $W$-invariant scalar product on $\Lambda$, which induces a metric  on $\La$. Using the $G$-action it can be continued to a metric $d$ on the whole of $X$ (see \cite{la}, 13.14 and \cite{brti1}, 7.4.20).

Note that if $n=2$, then $X$ is an infinite regular tree, with $q+1$ edges meeting in every vertex, where $q$ is the cardinality of the residue field. 

We denote by $X^0$ the set of vertices (i.e. $0$-dimensional faces) in $X$. We define a  simplex in $X^0$ to be a subset $\{x_1,\ldots, x_k\}$ of $X^0$ such that $x_1,\ldots,x_k$ are the vertices of a face in $X$. 

Let $\eta_i : \Gm  \rightarrow \bt$ be the cocharacter induced by 
mapping $x$ to the diagonal matrix with diagonal entries $d_1,\ldots,d_n$ such that $d_k= 1$ for $k \neq i$ and $d_i = x$. Then $\eta_1,\ldots, \eta_{n-1}$ is an $\R$-basis of $\Lambda$, and the set of vertices in $\La$ is equal to $\bigoplus_{i=1}^{n-1} \Z \eta_i$.

Let $\cL$ be the set of all homothety classes of $R$-lattices of full rank in $V$. We write $\{M\}$ for the class of a lattice $M$. Two different lattice classes $\{M'\}$ and $\{N'\}$ are called adjacent, if there are representatives $M$ and $N$ of $\{M'\}$ and $\{N'\}$ such that 
\[\pi N \subset M \subset N.\]
This relation defines a flag complex, namely  the simplicial complex with vertex set $\cL$ such that the simplices are the sets of pairwise adjacent lattice classes. We have a natural $G$-action on $\cL$ preserving the simplicial structure.

Moreover,  there is a $G$-equivariant bijection 
\[ \varphi: \cL \longrightarrow X^0\]
preserving the simplicial structures. 
If $\{N\} \in \cL$ can be written as $\{N\} = g\{M\}$ for some $g \in G$ and  $M = \pi^{k_1} R v_1 + \ldots + \pi^{k_n} R v_n$, then $\varphi(\{N\})$ is given by the pair $(g, \varphi\{M\}) \in G \times \La$, where 
\[\varphi(\{M\}) = \sum_{i=1}^{n-1}  (k_n-k_i) \eta_i\]
is a vertex in $\La$.

\section{Compactification of one appartment}
We write $\n$ for the set $\{1, \ldots, n\}$.
We continue to fix the base $v_1,\ldots, v_n$ of $V$ and the maximally split torus $\bt$ from section 2. Recall that $\Lambda = \bigoplus_{i = 1}^{n-1} \R \eta_i$, and that $\eta_n$ satisfies the relation $\eta_1 + \ldots + \eta_n = 0$. We will often write $\La = \sum_{i = 1}^n \R \eta_i$, bearing this relation in mind.

Let now $I$ be a non-empty subset of $\n$, and let $V_I$ be the subspace of $V$ generated by the $v_i$ for $i \in I$. We write $\bg^{V_I}$ for the subgroup of $\bg = PGL(V)$ consisting of the elements fixing the subspace $V_I$, and $\bg_I$ for the group $PGL(V_I)$. Then we have a natural restriction map
\[ \rho_I: \bg^{V_I} \longrightarrow \bg_I.\]
The torus $\bt$ is contained in  $\bg^{V_I}$, and its image under $\rho_I$ is a maximal $K$-split torus $\bt_I$ in $\bg_I$, namely the torus induced by the  diagonal matrices with respect to the base $\{v_i : i \in I\}$ of $V_I$. As usual, we write $T_I$, $G_I$  and $G^{V_I}$ for the groups of $K$-rational points. 

We put $\Lambda_I = X_\ast(\bt_I) \otimes_{\Z} \R$. Then $\rho_I$ induces a surjective homomorphism $\rho_{I \ast}: X_\ast(\bt) \rightarrow X_\ast (\bt_I)$, hence a surjective homomorphism of $\R$-vector spaces 
\[r_I : \Lambda \longrightarrow \Lambda_I.\]
For all $i \in I$ we write $\eta_i^I$ for the cocharacter of $\bt_I$ induced by mapping $x$ to the diagonal matrix with entry $x$ at the $i$-th place, i.e. $\eta_i^I = \rho_{I \ast} \eta_i$.

Then $\Lambda_I = \sum_{i \in I} \R \eta_i^I$, subject to the relation $\sum_{i \in I} \eta_i^I = 0$. In particular, $\La_{\{i\}} = 0$. Note that $r_I (\sum_{i = 1}^n x_i \eta_i) = \sum_{i \in I} x_i \eta_i^I$.

Let $\nu_I : T_I \rightarrow \Lambda_I$ be the unique homomorphism satisfying $< \nu_I(z), \chi> = -v(\chi(z))$ for all $\chi \in X^\ast( \bt_I)$. It is compatible with $\nu$, i.e. the following diagram is commutative:

$$
\begin{CD}
T @>{\nu}>> \Lambda \\
@V{\rho_I}VV @VV{r_I}V\\
T_I @>{\nu_I}>> \Lambda_I
\end{CD}
$$

Now we define
\[ \ov{\Lambda} = \Lambda \cup \bigcup_{\emptyset \neq I \subset \n}\Lambda_I = \bigcup_{\emptyset \neq I \subseteq \n} \Lambda_I.\]
Recall that we write ``$\subset$'' for a strict subset, and ``$\subseteq$'' if equality is permitted. Here $\La_{\n} = \La$ and $r_{\n}$ is the identity.

Let us now define a topology on $\Labar$. For all $I \subset \n$ we put
\[ D_I = \sum_{i \notin I} \R_{\geq 0} (-\eta_i).\]
We think of $D_I$ as a ``corner'' around $\Lambda_I$.
For all open and bounded subsets $U \subset \La$ we define
\[ C_U^I = (U + D_I) \cup \bigcup_{I \subseteq J \subset \n} r_J (U + D_I).\]
We take  as a base of our topology on $\Labar$ the open subsets of $\La$ together with these sets  $C_U^I$ for all non-empty $I \subset \n$ and all open bounded subsets $U$ of $\Lambda$. 

Note that every point $x \in \Labar$ has a countable fundamental system of neighbourhoods. This is clear for $x \in \La$. If $x$ is in $\La_I$ for some $I \subset \n$, then choose some $z \in \La$ with $r_I z = x$, and choose a countable decreasing fundamental system of bounded open neighbourhoods $(V_k)_{k \geq 1}$ of $z$ in $\La$. Put $U_k = V_k + \sum_{i \notin I} k (- \eta_i)$. This is an open neighbourhood of $z +  k \sum_{i \notin I} (- \eta_i)$. Then $(C_{U_k}^I)_{k \geq 1}$ is a fundamental system of open neighbourhoods of $x$. 

Hopefully the next result will shed some light on the definition of the topological space $\Labar$.

Recall from section 2, that we have a $G$-equivariant bijection $\varphi$ between equivalence classes of lattices of full rank in $V$ and vertices in the building $X$. If we restrict $\varphi$ to lattices which can be diagonalized with respect to $v_1,\ldots, v_n$, i.e. which have an $R$-basis consisting of multiples of these elements,  then we get a bijection between these diagonal lattices and vertices in $\La$. Applying this to the group $\bg_I = PGL(V_I)$, we get a bijection $\varphi_I$ between classes of diagonal lattices in $V_I$ with respect to the $v_i$ for $i \in I$, and vertices in $\La_I$.

\begin{prop} Let $(M_k)_{k \geq 1}$ be a sequence of diagonal lattices in $V$ and let $N$ be a diagonal lattice in $V_I$. 
The sequence of vertices $\varphi(\{M_k\})$ converges to the vertex $\varphi_I(\{N\}) \in \La_I$ in our topology on $\Labar$ iff after passing to a subsequence there are lattices $M_k'$ equivalent to $M_k$ such that $M_{k+1}' \subseteq M_k'$ and such that $\bigcap_{k} M_k'$ is equivalent to $N$. 
\end{prop}
{\bf Proof: }We can write $M'_k = \bigoplus_{i = 1}^n  \pi^{a_{i,k}} R v_i$ for some integers $a_{i,k}$. Since $M'_{k+1} \subset M'_k$, we have $a_{i,k+1} \geq a_{i,k}$. Therefore for all $i$ the sequence $a_{i,k}$ becomes stationary or goes to infinity, so that $\bigcap_{k} M'_k = \bigoplus_{i \in I'}  \pi^{a_i} R v_i$, where $I'$ is the set of all $i$, such that $a_{i,k}$ becomes stationary, i.e. $a_{i,k} = a_i$ for all $k$ big enough. Let us call this intersection module $N'$. It is a lattice in $V_{I'}$, and by assumption equivalent to $N$, so that $I = I'$.

Besides, we have $\varphi (\{M'_k\}) = \sum_{i = 1}^n (-a_{i,k}) \eta_i$, and  $\varphi_I(\{N'\}) = \sum_{i \in I} (-a_i )\eta_i^I$. If $k$ is big enough, we
have $\varphi(\{M'_k\}) = \sum_{i \in I} (-a_i) \eta_i + \sum_{i \notin I} (-a_{i,k}) \eta_i$, with $a_{i,k}$ arbitrarily large. If we take one of the systems of fundamental neighbourhoods of $\varphi_I(\{N\}) = \varphi_I(\{N'\})$ constructed previously, we find that every one of them must contain a point $\varphi(\{M'_k\})$, so that $\varphi(\{M_k\})$ converges indeed to $\varphi_I(\{N\})$.

To prove the other direction, assume that $\varphi(\{M_k\})$ converges to $\varphi_I(\{N\})$ for $M_k = \bigoplus_{i =1}^n \pi^{a_{i,k}} R v_i$ and $N = \bigoplus_{i \in I} \pi^{b_i} R v_i$. Looking at the fundamental neighbourhoods as above, we find that for any fixed $i_0 \in I$ the sequence $a_{i,k} - a_{i_0,k}$ is unbounded for $i \notin I$, and goes to $b_i - b_{i_0}$ for $i \in I$. This implies our claim.
 \hfill$\Box$

We will now show that the space $\Labar$ is compact.

Fix some $i \in \n$. 
We write $D_i$ for $D_{\{i\}} = \sum_{j \neq i} \R_{\geq 0} (-\eta_j)$,
the ``corner in $\La$ around the point $\La_{\{i\}}$''. Besides, let
\[ E_i ={D_i} \cup \bigcup_{i \in J \subset \n} r_J(D_i) \subset \Labar,\]
the ``closed corner in $\Labar$ around the point $\La_{\{i\}}$''.

\begin{lem}
i) Each point in $\Labar$ lies in one of the $E_i$.

ii) Each $E_i$ is  closed in $\Labar$.

\end{lem}
{\bf Proof: } i) Let $x = \sum_{j = 1}^n x_j (-\eta_j)$ be a point in $\La$.
Note that the relation $\sum_{j=1}^n \eta_j = 0$ implies that we can write $x = \sum_{j \neq i} (x_j - x_i) (-\eta_j)$ for all $i$. Now $E_i \cap \La = D_i$ is the set of all $x$ such that all $x_j-x_i$ are non-negative. In other words, a point $x = \sum_{j=1}^n x_j ( -\eta_j)$ is in $D_i$ iff $x_i$ is the minimum of all the coefficients $x_j$. This implies that for given $x$ we always find some $E_i$ containing it.
As similar argument holds if $x$ is contained in a boundary piece $\La_J$, since if $J$ contains $i$, the boundary piece $E_i \cap \La_J$ is the set of all $x = \sum_{j \in J} x_j (-\eta_j^J)$ such that $x_i$ is the minimum of all the $x_j$ for $j \in J$. (If $i$ is not contained in $J$, then of course $E_i \cap \La_J$ is empty.)

ii) Take some $x$ not contained in $E_i$. Then $x$ is in some $\La_J$ for $J \subseteq \n$ (possibly $\n$). Since the point $\La_{\{i\}}$ is contained in $E_i$, we know that $J \neq \{i\}$. Let us first assume that $i$ is contained in $J$. We write $x = \sum_{j \in J, j \neq i} x_j (-\eta_j^J)$. Since $x$ is not in $E_i$, our considerations in part i) imply that 
one of the $x_j$ for $j \in J$, say $x_{j_0}$, must be negative. The point $z = \sum_{j \in J, j \neq i}x_j (-\eta_j)$ in $\La $ projects to $x$, i.e. $r_J (z) = x$. This point must also be in the complement of $E_i$. Since $E_i \cap \La = D_i$ is closed, we find a bounded open neighbourhood $U$ of $z$ in $\La$ which is disjoint from $E_i$, and which contains only points $y = \sum_{j \neq i} y_j (-\eta_j)$ with $y_{j_0} < 0$. Now it is easy to check that the open neighbourhood $C_U^J$ of $x$ is also disjoint from $E_i$, which proves our claim. The remaining case $i \notin J$ can be treated with a similar argument.\hfill$\Box$

Let $R_{\geq 0, \infty}$ be the compactified half-line $\R_{\geq 0} \cup \{\infty\}$ with the topology generated by all intervals $[0,a[$, $]b,c[$ and $]b,\infty]$ for $a>0$ and $b,c \geq 0$. The space $\Rinfty$ is compact and contractible. A contraction map $r : \Rinfty \times [0,1] \rightarrow \Rinfty$ is given by 
\[r(x,t) = \frac{(1-t)x}{1+tx}\mbox{ for } x \in \R, \quad \mbox{and } \quad r(\infty,t)= \left\{\begin{array}{ll}\infty, &\mbox{if } t = 0 \\ \frac{1-t}{t}, & \mbox{if } t \neq 0 \end{array}\right. ,\]
see \cite{la}, 2.1.
Let us fix some $i \in \n$. We will now compare $E_i$ to $\Rinfty^{n-1}$, which we write as $\bigoplus_{j \neq i} \Rinfty e_j$ for a basis $e_j$.
Recall from the proof of Lemma 3.2 that in the case $i \in I$ we can describe $E_i \cap \La_I$ as the set of all $x = \sum_{j \in I, j \neq i} x_j (-\eta_j^I)$ with non-negative $x_j$. Hence the following map is a bijection
\begin{eqnarray*}
\alpha_I : & E_i \cap \Lambda_I & \rightarrow \{\sum_{j \neq i} x_j e_j \in \bigoplus_{j \neq i} \Rinfty e_j: x_j = \infty \mbox{ iff }j \notin I\}\\
~ & \sum_{j \in I, j \neq i} x_j (-\eta_j^I) & \mapsto  \sum_{j \in I, j \neq i} x_j e_j + \sum_{j \notin I} \infty\,  e_j.
\end{eqnarray*}
For $I = \n$ the map $\alpha_{\n}: E_i \cap \La = D_i \rightarrow \bigoplus_{j \neq i} \R_{\geq 0} e_j$ can be continued to a homomorphism of $\R$-vector spaces
\[ \alpha_\La: \La \longrightarrow \bigoplus_{j \neq i} \R e_j,\]
which is a homeomorphism. Putting all the maps $\alpha_I$ together, we get a bijection
\[ \alpha: E_i \longrightarrow  \bigoplus_{j \neq i} \Rinfty e_j,\]
whose restriction to $E_i \cap \Lambda$ is a homeomorphism. We even have the following fact: 

\begin{lem}
With respect to the topology on $E_i$ induced by $\Labar$, the map $\alpha$ is a homeomorphism on the whole of $E_i$.
\end{lem}

{\bf Proof: }For all $j \neq i$ choose an open interval $A_j$ in $\Rinfty$, which is either of the form $A_j = [0,a_j[$ or $A_j = ]b_j,c_j[$ or of the form $A_j= ]b_j,\infty]$. We claim that the preimage of $A = \sum_{j \neq i} A_j e_j$ is open in $E_i$. 

We put $A_j' = ]-1,a_j[$, if $A_j = [0,a_j[$. In all the other cases we put $A_j' = A_j$. Let $A' = \sum_{j \neq i} A'_j e_j$ and put
\[ W = \alpha_\La\inv(A' \cap \bigoplus_{j \neq i} \R e_j).\]
Since $\alpha_\La$ is a homeomorphism, $W$ is open in $\Lambda$. Obviously, we have $W \cap E_i = \alpha\inv (A) \cap \Lambda$. Now put 
\[ I = \{ j \in \n : \infty \notin A_j\} \cup \{i\} .\]
We can assume that $I \neq \n$.
Choose some positive real number $b$ strictly bigger than all the $b_j$ for $j \notin I$. Then $U = W \cap \{ x= \sum_{j \neq i} x_j (-\eta_j): x_j < b \mbox{ for }j \notin I\}$ is an open bounded subset of $\Lambda$. Note that $U + D_I = W$. We claim that $\alpha\inv(A) = C_U^I \cap E_i$.

Indeed, every element $u$ in $W = U+D_I$ can be written as $u = \sum_{j \neq i} x_j (-\eta_j)$ with  $x_j \in A'_j$. Let $J$ be a subset of $\n$ containing $I$. If $r_J (u) = \sum_{j \in J, j\neq i}x_j (-\eta_j^J)$ is in $E_i$, we have
$x_j \in A_j$ for all $j \in J$ not equal to $i$. This implies that $\alpha(r_J(u))$ is contained in $A$. On the other hand, let $y = \sum_{j \neq i} y_j e_j$ be an element in $A$, i.e. $y_j \in A_j$. 
Put $J = \{j \neq i: y_j \neq \infty\} \cup \{i\}$. Then $J$ contains $I$. We put $x_j = y_j$, if $j \neq i$ is in $J$. If $j \notin J$, we choose an arbitrary element in $A_j \cap \R$ and call it $x_j$. Then $x = \sum_{j \neq i} x_j (-\eta_j)$ is contained in $W \cap D_i$, so that $r_J(x)$ is an element in $C_U^I\cap E_i$ which satisfies $\alpha(r_J(x)) = y$. Hence we also find $\alpha\inv (A) \subset C_U^I \cap E_i$.

Therefore $\alpha$ is continuous. It remains to show that $\alpha$ is open. 
Let $U$ be an open, bounded subset of $\Lambda$ and $I \subset \n$ non-empty. We will show that $\alpha(C_U^I \cap E_i)$ is open. Let $x$ be a point in $C_U^I \cap E_i$ lying in $\Lambda_J$ for some J containing $I$ and $i$. Hence $x = \sum_{j \in J, j \neq i} x_j (- \eta_j^J)$ with non-negative $x_j$. Since $x$ is contained in $C_U^I$, we can find some $z = \sum_{j \neq i} z_j (-\eta_j)$ in $U + D_I$ such that $r_J(z) = x$ (hence $z_j \geq 0$ for $j \in J$) and $z_j >0$ for $j \notin J$. Then $z \in (U+ D_I) \cap E_i$. 

Since the restriction of $\alpha$ to $E_i \cap \La$ is open, we find open intervals $A_j$ (of the form $[0,a_j[$ or $]b_j,c_j[$) in $\R_{\geq 0 }$ such that $A = \sum_{j \neq i} A_j e_j$ is an open neighbourhood of $\alpha(z)$ in $\bigoplus_{j \neq i} \R_{\geq 0} e_j$ which is contained in $\alpha((U+D_I) \cap E_i)$. We can also assume that for $j \notin J$ the interval $A_j$ does not contain $0$.

Now put $A'_j = A_j$ if $j\neq i$ is contained in $J$, and put $A'_j = ]b_j, \infty]$ if $j$ is not contained in $J$, and $A_j = ]b_j,c_j[$. (The interval $A_j$ looks indeed like this since we took care to stay away from zero.)

It is easy to see that $A' = \sum_{j \neq i} A_j' e_j$ is contained in  $\alpha (C_U^I \cap E_i)$.
Hence we found an open neighbourhood $A'$ around $\alpha(x)$ in $\alpha (C_U^I \cap E_i)$.\hfill$\Box$

\begin{thm}
The topological space $\Labar$ is compact and contractible, and $\La$ is an open, dense subset of $\Labar$. 
\end{thm}
{\bf Proof:} By the previous result, all $E_i$ are compact and contractible. Since $\Labar$ is the union of the $E_i$, it is also compact. One can check that the contraction maps are compatible, so that $\Labar$ is contractible. It is clear that $\La$ is open and dense in $\Labar$.\hfill$\Box$

Our next goal is to extend the action of $N$ on $\La$ to a continous action on the compactification $\Labar$. Recall that we identified $W$ with the group of permutation matrices in $N$, so that  $N = T \rtimes W$. For $w \in W$ we denote the induced permutation of the set $\n$ also by $w$, i.e. we abuse notation so that $w(v_i) = v_{w(i)}$. 

Let $I$ be a non-empty subset of $\n$. We define a map 
\[w: \La_I \longrightarrow \La_{w(I)}\]
by sending $\eta_i^I$ to $\eta_{w(i)}^{w(I)}$. This gives an action of $W$ on $\Labar$. Note that it is compatible with $r_J$, i.e. we have
\[ w \circ r_J = r_{w(J)} \circ w\]
on $\La$. Besides, we can combine the maps $\nu_I: T_I \rightarrow \Lambda_I$ with the restriction map $\rho_I : \bt \rightarrow \bt_I$ to define a map
$\nu_I \circ \rho_I: T \rightarrow \Lambda_I$, so that $T$ acts by affine transformations on $\Lambda_I$. Recall  that $r_I(\nu(t)) = \nu_I (\rho_I (t))$ for all $t \in T$.

It is easy to check that these two actions give rise to an action of $N = T \rtimes W$ on $\Labar$, which we denote by $\nu$. 
\begin{lem}
The action $\nu: N \times \Labar \longrightarrow \Labar$ is continuous and extends the action of $N$ on $\La$.
\end{lem}
{\bf Proof: }Let first $w$ be an element of $W$, and let $C_U^I$ one of our open basis sets. Then $\nu(w)(U + D_I) = \nu(w) (U) + D_{w(I)}$, since $\nu(w)$ is a linear map on $\La$. Besides, we have $\nu(w) (r_J (U + D_I)) = r_{w(J)} (\nu(w)(U) + D_{w(I)})$, so that $\nu(w)(C_U^I) = C_{\nu(w)(U)}^{w(I)}$.

Now take some element $t \in T$. Then $\nu(t) (U + D_I) = \nu(t)(U) + D_I$, since $\nu(t)$ acts by translation. Besides, we have $\nu(t)(r_J(U + D_I)) = r_J(\nu(t) (U+D_I))$, so that $\nu(t)(C_U^I) = C_{\nu(t)(U)}^I$.

Hence for all $n \in N$ the action $\nu(n)$ on $\Labar$ is continuous. Since the kernel of the map $\nu: T \rightarrow \La$  is an open  subgroup of $N$ (see \cite{la}, Prop. 1.2), which obviously acts trivially on $\Labar$, we find that the action is indeed continuous.\hfill$\Box$
 
\section{Compactification of the building}
We can now define the compactification of the building $X$. For all non-empty subsets $\Omega$ of $\Labar$ and all roots $a \in \Phi$ we put
\begin{eqnarray*}
 f_\Omega(a) & =  \inf\{t: \Omega \subseteq\overline{\{ z \in \La: a(z) \geq -t\}} \} \\
~ & = -\sup \{ t: \Omega \subseteq \ov{\{ z \in \La: a(z) \geq t \} } \}
\end{eqnarray*}
Here we put $\inf \emptyset = \sup \R = \infty$ and $\inf \R = \sup \emptyset = - \infty$. Moreover, if $\Omega = \{x\}$ consists of one point only, then we write $f_x(a) = f_{\{x\}}(a)$. Note that 
\[ f_x(a) = -a(x) \quad \mbox{ for all }x \in \La,\]
and that 
\[f_{\om_1}(a) \leq f_{\om_2}(a), \quad \mbox{if }\om_1 \subseteq \om_2.\] 
Recall our generalized valuation map $\psi_a: U_a \rightarrow \Z \cup \{\infty\}$ from section 2. We can now define a subgroup 
\[ U_{a,\Omega} = U_{a,f_\Omega(a)} = \{ u \in U_a : \psi_a(u) \geq f_\Omega(a)\}\]
of $U_a$, where $U_{a, \infty} = 1$ and $U_{a, -\infty} = U_a$. By $U_\Omega$ we denote the subgroup of $G$ generated by all the $U_{a,\Omega}$ for roots $a \in \Phi$. Note that if $\Omega= \{x\}$ for some point $x \in \La$, then this coincides with our previous definition of $U_x$. 
We will now investigate these groups $U_x$ for boundary points of $\Labar$.

\begin{prop}
Recall that we denote by $a_{ij}$ the  root of $T$ induced by the character $\chi_i - \chi_j$.  Put $a = a_{ij}$, and let $x$ be a point in $\La_I$ for some $I \subseteq  \n$. 

i) If $j \notin  I$, we have $f_x(a) = - \infty$, so that $U_{a,x}$ is the whole group $U_a$. 

ii) If $j \in I$ and $i \notin I$, then we have $f_x(a) = \infty$, so that $U_{a,x} = 1$. 

iii) If $i$ and $j$ are contained in $I$, then $a$ is equal to $\rho_I^\ast (b)$ for some root $b$ of the torus $T_I$ in $G_I$. In this case we have
$f_x(a) = -b(x)$. For any $z \in \La$ with $r_I(z) = x$ we also have $f_x(a) = -a(z)$. 
\end{prop}
{\bf Proof: } i) Choose some $z \in \La$ such that $r_I(z) =x$. If $i$ is contained in $I$, put $z_k = z + \sum_{l \notin I} k (-\eta_l)$. If $i$ is not in $I$, then we define $z_k = z + \sum_{l \notin I, l \neq j} k (-\eta_l) - 2k \eta_j$. In both cases we find that $a(z_k)$ equals $a(z) + k$, hence it goes to infinity. Since the $z_k$ converge to $x$, we find that $x$ lies indeed in the closure of any set of the form $\{a \geq s\}$, which implies our claim.

ii) We choose again some $z \in \La$ with $r_I(z) = x$. Let $V_k$ be a countable decreasing fundamental system of bounded open neighbourhoods of $z$. This defines a fundamental system of open neighbourhoods $C_{U_k}^I$ around $x$, where $U_k = V_k + \sum_{i \notin I} k (-\eta_i)$. Now suppose that $x$ is contained in the closure of the set $\{z \in \Lambda: a(z) \geq s\}$. Then we find for all $k$ some $y_k$ in $C_{U_k}^I \cap \Lambda$ satisfying $a(y_k) \geq s$. We can write $y_k = z_k + \lambda_k$ for some $z_k \in V_k$ and $\lambda_k = \sum_{l \notin I } \lambda_{k,l} (-\eta_l)$ with $\lambda_{k,l}  \geq k$. Now $a(z_k)$ is bounded, but $a(\lambda_k) = -\lambda_{k,i}$, so that $a(y_k)$ cannot be bounded from below. Hence we find indeed that $f_x(a)$ must be $\infty$.

iii) Recall that $\bg_I$ is the group $PGL(V_I)$, and $\bt_I$ is the maximal $K$-split torus induced by the diagonal matrices with respect to the $v_i$ for $i \in I$. Then the root system  corresponding to $\bt_I$ and $\bg_I$ is 
\[\Phi_I = \{ b_{ij}: i \neq j \mbox{ in } I\}\]
where $b_{ij}$ is the character mapping a diagonal matrix with entries $t_i$ for $i \in I$ to $t_i / t_j$. 
Hence it is clear that in our case $i,j \in I$ the root $a = a_{i,j}$ of $T$ is induced by the root $b = b_{ij}$ of $T_I$. Note that this implies that for all $z \in \La$ we have $a(z) = b(r_I(z))$.

It suffices to show that 
\[ \ov{\{ z \in \La: a(z) \geq s\}} \cap \La_I = \{ x \in \La_I: b(x) \geq s\}.\]
Take some $x$ contained in the left hand side, and choose some $z \in \La$ with $r_I(z) = x$. Besides, we take again fundamental neighbourhoods $V_k$ around $z$ and use them to construct the open neighbourboods $C_{U_k}^I$ around $x$. Each $C_{U_k}^I$ must contain some $y_k \in \La$ satisfying $a(y_k) \geq s$. Note that we can write $y_k = z_k + \lambda_k$, where $z_k$ is in $V_k$ and $\lambda_k$ is a linear combination of $\eta_i$ for $i \notin I$. Besides $a(y_k) = a(z_k)$, so that the sequence of $a(y_k)$ converges to $a(z)$. Since all $a(y_k)$ are $\geq s$, we find $b(x) = b(r_I z) = a(z) \geq s$. 

On the other hand, suppose that $x$ is a point in $\La_I$ satisfying $b(x) \geq s$. Again, we choose some $z \in \La$ with $r_I(z) = x$, and neighbourhoods $V_k$ around $z$. For any $k$ the point
$z_k = z + k \sum_{ l \notin I}( - \eta_l)$ lies in $C_{U_k}^I$. Besides, $a(z_k) = a(z) = b(r_I(z)) = b(x)$ is bounded below by $s$. This implies that $x$ lies indeed in the closure of $\{ z \in \La: a(z) \geq s\}$.\hfill$\Box$

\begin{prop} Let $x$ be in $\La_I$ and let $a= a_{ij} \in \Phi$ be a root.

i) Each $U_{a,x}$ (and hence also $U_x$) leaves the vector space $V_I$ invariant. Hence $U_{a,x}$ is contained in $G^{V_I}$.

ii) If $i$ and $j$ are not both in $I$, we have $\rho_I(U_{a,x}) = 1$. If $i$ and $j$ are both in $I$, and the root $a$ is induced by the root $b$ of $T_I$, then $\rho_I$ induces an isomorphism $U_{a,x} \rightarrow   U_{b,x}^I$, where $U_{b,x}^I$ is defined with the root group $U_b^I$ in $\bg_I$ as described in section 2.
\end{prop}

{\bf Proof: } Recall that $u \in U_a$ maps $v_l$ to itself, if $l$ is not equal to $j$, and it maps $v_j$ to $v_j + \omega v_i$, where $\psi_a(u) = v(\omega)$. Hence our claim in i) is clear if $j$ is not contained in $I$ or if both $i$ and $j$ are contained in $I$. In the remaining case we saw in 4.1 that $U_{a,x}$ is trivial, so that i) holds in any case.

Let us now prove ii). If $j$ is not contained in $I$, then each $u \in U_a$ induces the identity map on $V_I$. If $j$ is in $I$, but $i$ is not, then $U_{a,x}$ is trivial. Hence in both cases we find that $\rho_I(U_{a,x}) = 1$. 
Let us assume that both $i$ and $j$ are contained in $I$, and let $u$ be an element of $U_{a,x}$. Then $\rho_I(u) \in PGL(V_I)$ is induced by the matrix mapping $v_l$ to $v_l$ for all $l \neq j$ in $I$, and $v_j$ to $v_j + \omega v_i$ with some $\omega$ having valuation $\geq f_x(a)$. By our description of the groups $U_b^I$ in section 2 we find that $\rho_I(u)$ is contained in $U_b^I$ and has valuation $\psi_b(\rho_I x) = v(\omega) \geq f_x(a)$. By 4.1, $f_x(a) = -b(x)$, so that $\rho_I(u)$ lies indeed in $U_{b,x}^I$. The homomorphism $\rho_I: U_{a,x} \rightarrow U_{b,x}^I$ is obviously bijective.\hfill$\Box$

Note that the map $x \mapsto f_x(a)$ is in general not continuous on $\Labar$. Take some $z \in \La$ and define a sequence $x_k = z + \sum_{i \notin I} k (-\eta_i)$ for some non-empty $I$ such that the complement $\n \backslash I$ contains at least two elements $i$ and $j$. Then $a=a_{ij}$ has the property that $a(x_k) = a(z)$, so that $f_{x_k}(a)$ is constant. But the sequence $x_k$ converges to the point $r_I(z)$ in $\La_I$, for which $f_{r_I(z)}(a) = - \infty$ holds. 

Nevertheless, we have the following result:
\begin{lem}
Let $x_k$ be a sequence of points in $\Labar$, which converges to $x \in \Labar$. Let $u_k \in U_{a,x_k}$ be a sequence of elements, converging to some $u$ in the big group $U_a$. Then $u$ lies in fact in $U_{a,x}$.
\end{lem}
{\bf Proof: }Note first of all, that the statement is clear if $f_x(a) = - \infty$, since then $U_{a,x} = U_a$. It is also clear if $f_{x_k}(a)$ converges to $f_x(a)$, since the map $\psi_a : U_a \rightarrow \Z \cup \{\infty\}$ is continuous. 
Assume that $f_x(a) = \infty$. Then any set $\{\ov{z: a(z) \geq s}\}$ 
 contains only finitely elements $x_k$. This implies that the sequence $f_{x_k}(a)$ goes to $\infty = f_x(a)$, so that in this case our claim holds by continuity.

The only case which is left is that $f_x(a)$ is real. Assume that $x \in \La_I$. By 4.1, we must have $a=a_{ij}$ with $i$ and $j$ in $I$. Choose some $z \in \La$ with $r_I(z) = x$, and a decreasing fundamental system of open neighbourhoods $V_k$ around $z$. As before, we use them to define neighbourhoods $C_{U_k}^I$ around $x$. Since $x_k$ converges to $x$, we can assume that $x_k$ is contained in $C_{U_k}^I$. Then $x_k \in \La_{J_k}$ for some $J_k$ containing $I$. By definition of $C_{U_k}^I$ we find some $y_k \in V_k$ and some coefficients $\alpha_l$ such that $z_k = y_k + \sum_{l \notin I} \alpha_l (-\eta_l)$ satisfies $r_{J_k} (z_k ) = x_k$. By  4.1  we have $f_{x_k}(a) = -a(z_k) = -a(y_k)$ and $f_x(a) = -a(z)$. Hence $f_{x_k}(a)$ converges to $f_x(a)$, and our claim follows again by continuity.\hfill$\Box$

\begin{prop}
For $x \in \Labar, n \in N$ and $a \in \Phi$ we have
\[n U_{a,x} n\inv = U_{\ov{n}(a), \nu(n)(x)} ,\]
where $\nu$ denotes the action of $N$ on $\Labar$, and $n \mapsto \ov{n}$ denotes the quotient map from $N$ to the Weyl group $W$ (which acts on the roots).
In particular, we have $n U_x n\inv = U_{\nu(n)(x)}$.
\end{prop}
{\bf Proof: }Fix some $n \in N$ and denote by $p$ the permutation matrix mapping to $\ov{n}$, i.e. $n = t p$ for some $t \in T$. We denote by $p$ also the corresponding permutation of $\n$. If $a = a_{ij}$, then $\ov{n}(a) = a_{p(i) p(j)}$. By 4.1  we find that if $x \in \La_I$ and  $j \notin I$ both $f_x(a)$ and $f_{\nu(n)(x)}(\ov{n} a)$ are equal to $- \infty$. Since $n U_a n\inv = U_{\ov{n} a} $, our claim holds in this case. Similarly, if $j \in I$ and $i \notin I$ both $f_x(a)$ and $f_{\nu(n)(x)}(\ov{n} a)$ are equal to $\infty$, so that both $U_{a,x}$ and $U_{\ov{n}a, \nu(n)(x)}$ are trivial.

We can therefore assume that $x$ is in some $\La_I$ such that both $i$ and $j$ are contained in $I$. Recall that we composed the action of $N$ on $\La$ from the natural action of the Weyl group on $\La$ and translation by $\nu(t)$ for elements in the torus $T$. Hence for all $z \in \La$ we have $\nu(n)(z) = \nu(n)(0) + \ov{n}(z)$. We can now calculate
\begin{eqnarray*}
a(\nu(n)(z)) & = & a(\nu(n)(0)) + a(\ov{n}(z)) \\
~ & = & a(\nu(n)(0)) + (\ov{n}\inv a)(z), \end{eqnarray*}
so that 
$\nu(n\inv)\{ z\in \La: (\ov{n}a) (z) \geq s\} = \{z \in \La: a(z) \geq a(\nu(n\inv)(0)) +s\}$.
Since $\nu(n\inv)$ is a homeomorphism, we also have 
\[ \nu(n\inv)\ov{\{ z: \ov{n} a(z) \geq s\}} = \ov{\{z: a(z) \geq a(\nu(n\inv)(0)) +s\}}.\]
Hence 
\begin{eqnarray*}
f_{\nu(n)(x)}(\ov{n}(a)) & = & \inf\{t: \nu(n)(x) \in \ov{\{ \ov{n}(a) \geq -t\} }\}\\
~& = & \inf\{ t: x \in \ov{\{z: a(z) \geq a(\nu(n\inv)(0)) -t\} }\}\\
~ & = & f_x(a) + a(\nu(n\inv)(0)).
\end{eqnarray*}
Since we have $n U_{a,s} n\inv = U_{\ov{n}(a), s + a(\nu(n\inv)(0))}$ for all real numbers $s$ (see \cite{la}, 11.6), we find that indeed
\[ n U_{a,x}n\inv = U_{\ov{n}(a), f_x(a)+ a(\nu(n\inv)(0))} = U_{\ov{n}(a), f_{\nu(n)(x)}(\ov{n}(a))}= U_{\ov{n}(a), \nu(n)(x)},\]
whence our claim.\hfill$\Box$

Recall that - forgetting about the special nature of our ground field - our root system $\Phi$ in $\La^\ast$ defines a finite set of hyperplanes in $\La^\ast$ and therefore a decomposition of $\La^\ast$ into faces (see \cite{bou},V,1). The maximal faces are called (spherical) chambers. Any chamber defines an order on $\La^\ast$ (\cite{bou}, VI, 1.6). We denote the positive roots with respect to this order by $\Phi^+ = \Phi^+(C)$, and the negative roots by $\Phi^-= \Phi^-(C)$. In fact, for any subset $\Psi$ of $\Phi$ such that $\Psi$ is additively closed and $\Phi$ is the disjoint union of $\Psi$ and $-\Psi$, there exists a chamber $C$ such that $\Psi= \Phi^+(C)$ (see \cite{bou}, VI, 1.7). In particular, $\Phi^-$ is the set of positive roots for a suitable chamber. 

The following lemma will be useful to reduce claims about the groups $U_{a,\om}$ for $\om \subset \Labar$ to claims about subsets of $\La$, where we can apply the ``usual'' theory of the Bruhat-Tits building.
\begin{lem}
Let $\om \subset \Labar$ be a non-empty set, and fix some chamber $C$. Put $\Phi^+ = \Phi^+(C)$. Assume that for every $a \in \Phi^+$ we have $m$ elements
\[u_{a,1}, \ldots, a_{a,m} \in U_{a,\om},\]
such that at least one of all these $u_{a,i}$'s is non-trivial. 
Then there exists a non-empty subset $\Omega'$ of $\La$ such that $u_{a,i} \in U_{a,\om'}$ for all $i = 1,\ldots, m$ and such that $\om \subset \ov{\om}'$.
In particular, we have $U_{a,\om'} \subset U_{a,\om}$ for all roots $a \in \Phi$.
\end{lem}
{\bf Proof: }We denote by $l_a$ the infimum of all $\psi_a(u_{a,i})$ for $i = 1,\ldots,m$. If all $u_{a,i}$ are trivial, then $l_a = \infty$. This cannot happen for all $a \in \Phi^+$. 
Then we put 
\[\om_a' = \{\ov{z \in \La: a(z) \geq -l_a}\}.\]
(If $l_a = \infty$, then $\om_a' = \Labar$.) Besides, put $\om' = \La \cap \bigcap_{a \in \Phi^+}\om_a'$. Note that $\om'$ contains the intersection of all sets 
$\{ z \in \La: a(z) \geq -l_a\}$ for $a \in \Phi^+$. If $l$ is the minimum of all the $l_a$, then this set contains all $z \in \La$ satisfying $a(z) \geq -l$ for all $a \in \Phi^+$. By looking at a base of $\Phi$ corresponding to $\Phi^+$, we see that such $z$'s exist. Hence $\om'$ is non-empty. 

By construction, $f_{\om_a'}(a) = l_a$, and the inclusion $\om' \subset \om'_a$ gives us $f_{\om'}(a) \leq f_{\om'_a}(a)$. Therefore $\psi_a(u_{a,i}) \geq l_a \geq f_{\om'}(a)$, which implies that $u_{a,i}$ is indeed contained in $\Uaomp$ for all $i = 1,\ldots,m$. 
It remains to show that $\om \subset \ov{\om}'$. Since $f_{\om'}(a) = f_{\ov{\om}'}(a)$ for all roots $a$, this implies that $f_\om (a) \leq f_{\om'}(a)$, hence $\Uaomp \subset \Uaom$ for all $a \in \Phi$.  

We are done if we prove the following claim:

\[(\ast) \mbox{ For any } \Psi \subseteq \Phi^+ \mbox{ and real numbers } s_a \mbox{ we have }\bigcap_{a \in \Psi}\{ \ov{a \geq s_a}\} = \ov{\bigcap_{a \in \Psi} \{ a \geq s_a\} }.\]

It is clear that the right hand side is contained in the left hand side. So suppose $x$ is an element in $\bigcap_{a \in \Psi}\{ \ov{a \geq s_a}\}$. Let $I$ be the subset of $\n$ such that $x \in \La_I$. We choose a system of open neighbourhoods $V_k$ of some point  in $\La$ projecting to  $x$ and construct $C_{U_k}^I$. We are done if we can show that any $C_{U_k}^I$ intersects  $\bigcap_{a \in \Psi} \{ a \geq s_a\}$ non-trivially. Let $z_k$ be a point in $U_k = V_k + k \sum_{l \notin I} (-\eta_l)$ with $r_I(z_k) = x$. Besides, let $s_k$ be the maximum of $0$ and all the numbers $s_a-a(z_k)$ for all $a \in \Psi$. 

Note that $\Phi^+$ defines a linear ordering of the set $\n = \{1,\ldots, n\}$, namely $i \prec j$, iff $a_{ij} \in \Phi^+$. Hence there is a permutation $\pi$ of $\n$ satisfying $\pi(1) \prec \pi(2) \prec \ldots \prec \pi(n)$. Put $z = z_k - \sum_{l \notin I} (k + \pi\inv(l) s_k) \eta_l$. This is an element of $C_{U_k}^I$. It remains to  show that indeed $a(z) \geq s_a$ for all $a \in \Psi$. 

Let $a= a_{ij}$ be a root in $\Psi$. If both $i$ and $j$ are in $I$, we can apply 4.1 to deduce  $a(z) = a(z_k) = -f_x(a) \geq s_a$. Since $x$ is contained in $\{ \ov{a \geq s_a}\}$, it cannot happen that $j$ is in $I$, but $i$ is not. If $j$ is not in $I$, we find that
\[a(z) = \left\{ \begin{array}{rll}
a(z_k) + k  + \pi\inv(j) s_k & \geq s_a & \mbox{  ,if } i \in I \\
a(z_k) + (\pi\inv(j) - \pi\inv(i)) s_k & \geq s_a& \mbox{  ,if } i \notin I \end{array}, \right.
\]
since $a \in \Phi^+$ implies that $\pi\inv(i) < \pi\inv(j)$. Hence we get $a(z) \geq s_a$ for all $a \in \Psi$, which proves $(\ast)$.\hfill$\Box$

\begin{cor}
Assume that $a$ and $b$ are roots in $\Phi$ which are not linear equivalent (i.e. $a \neq \pm b$), and so that $a+b$ is in $\Phi$.  
If both $f_\om(a)$ and $f_\om(b)$ are  real numbers, then
\[f_\om(a+b) \leq f_\om(a) + f_\om(b).\]
If $f_\om(a) = -\infty$ and $f_\om(b) \neq \infty$, then $f_\om(a+b) = -\infty$.
\end{cor}
{\bf Proof: }By $(\ast)$ in the proof of the preceeding lemma, we have
\[\{\ov{a \geq s}\} \cap \{ \ov{b \geq r}\} = \{\ov{a \geq s, b \geq r}\} \subset \{\ov{a+b\geq s+r}\},\]
which implies our claim.\hfill$\Box$

Now we prove a statement about the structure of the groups $U_\om$ which will we crucial for our later results.

For $\Phi^+ = \Phi^+(C)$  we denote by ${\bf{U}}_{\Phi^+}$ the corresponding subgroup of $\bg$ (see \cite{bo}, 21.9), and by $U_{\Phi^+}$ the set of $K$-rational points. Similarly, we have ${\bf{U}}_{\Phi^-}$ and  $U_{\Phi^-}$.
For any non-empty subset $\om$ of $\Labar$ define
\[ U_\om^+ = U_{\Phi^+} \cap U_\om \quad \mbox{and} \quad U_\om^- = U_{\Phi^-} \cap U_\om.\]
Of course, these groups depend on the choice of some chamber $C$. We can use them to get some information about $U_\om$. 
\begin{thm}
i) The multiplication map induces a bijection $\prod_{a \in \Phi^\pm} \Uaom \longrightarrow U_\om^\pm$, where the product on the left hand side may be taken in arbitrary order.

ii) $U_a \cap U_\om = \Uaom$ for all $a \in \Phi$.

iii) $U_\om = U_\om^- U_\om^+ (N \cap U_\om)$.
\end{thm}

{\bf Proof: }Note that by \cite{la}, 12.5, our claim holds for all non-empty subsets $\om \subset \La$. Now take $\Omega \subseteq \Labar$, and denote by $L_a$ the group generated by $U_{a,\om}$ and $U_{-a, \om}$, and by $Y$ the subgroup of $N$ generated by all $N \cap L_a$ for $a \in \Phi$. For all $a \in \Phi^+$ choose an element $u_a \in \Uaom$. By 4.5 we find a subset $\om'$ of $\La$ such that $u_a \in \Uaomp$. Hence by \cite{la}, 12.5, the product of the $u_a$ in arbitrary order lies in $U_{\om'}^+ \subset U_\om^+$.

A similar argument using 4.5 shows that the image of $\prod_{a \in \Phi^+}U_{a,\om}$ under the multiplication map is indeed a subgroup of $U^+_\om$, which is independent of the ordering of the factors. We denote it by $H^+$. Similarly, we define the subgroup $H^-$ of $U_\om^-$ as the image of $\prod_{a \in \Phi^-} U_{a,\om}$ under the multiplication map. Now we can imitate the argument in \cite{la}, Proposition 8.9, to prove that the set $H^- H^+ Y$ does not depend on the choice of the chamber defining $\Phi^+$ and is invariant under multiplication from the left by $Y$ and $U_{a,\om}$ for arbitrary roots $a \in \Phi$. Hence $H^- H^+ Y = U_\om$. 

Since $U_{\Phi^+} \cap U_{\Phi^-} = \{1\}$ and $N \cap U_{\Phi^+} U_{\Phi^-} = \{1\}$ (by \cite{boti}, 5.15), we find $U_\om^- = (H^- H^+ Y) \cap U_{\Phi^-} = H^-$ and $U_\om^+ = H^+$, which proves i) and ii). Similarly, $N \cap U_\om = Y$, whence iii).\hfill$\Box$

For any subset $\om$ of $\Labar$ we write $N_\om = \{n \in N: \nu(n) x = x \mbox{ for all } x \in \om\}$. Besides, put
\[ P_\om= U_\om N_\om = N_\om U_\om,\]
which is a group since as in  4.4 one can show that $N_\om$ normalizes $U_\om$. If $\om = \{x\}$ we write $P_\om = P_x$.
 
We can now also describe the groups $P_\om$ for any non-empty subset $\Omega$ of $\Labar$:

\begin{cor} Fix some  $\Phi^+ = \Phi^+(C)$ as above, and let $\Omega$ be a non-empty subset of $\Labar$.

i) $P_\om = U_\om^- U_\om^+ N_\om = N_\om U_\om^+ U_\om^-$.

ii) $P_\om \cap U_{\Phi^\pm} = U_\om^\pm$ and $P_\om \cap N = N_\om$

\end{cor}

{\bf Proof: } i) The first equality follows from part iii) of the Theorem, if we show that $N \cap U_\om \subset N_\om$. It suffices to show for each root $a$ that $N \cap L_a \subset N_\om$. If both $f_\om(a)$ and $f_\om(-a)$ are real numbers, this follows from \cite{la}, 12.1. If $f_\om(a) = \infty$ of if $f_\om(-a) = \infty$, then our claim is trivial. Note that if $f_\om(a) = -\infty$, then $f_\om(-a)$ is either $\infty$ (then we are done) or $- \infty$. Hence the only remaining case is $f_\om(a) = f_\om(-a) = - \infty$. If $a = a_{ij}$, 4.1 implies that $\om \cap \La_J$ can then only be non-empty if $i$ and $j$ are not contained in $J$. Now 4.2 implies our claim. 

ii) Let us show first that $P_\om \cap U_{\Phi^-} = U^-_\om$. Obviously, the right hand side is contained in the left hand side. Take some $u \in P_\om \cap U_{\Phi^-}$. Using i), we can write it as $u^- u^+ n$ for $u^\pm \in U_\om^\pm$ and $n \in N_\om$. Then $n$ must be in $U_{\Phi^+} U_{\Phi^-} \cap N$, which is trivial by \cite{boti}, 5.15. Hence $u^+$ is contained in $U_{\Phi^-} \cap U_{\Phi^+}$, which is also trivial. Therefore $u= u^- \in U_\om^-$. The corresponding statement for the $+$-groups follows by taking $\Phi^-$ as the set of positive roots. It remains to show $P_\om \cap N = N_\om$. Take $u \in P_\om \cap N$. Then we write it again as $u = u^- u^+ n$. Hence $u^- u^+$ is contained in $U_{\Phi^-} U_{\Phi^+} \cap N$, which is trivial, so that $u = n \in N_\om$, as claimed.\hfill$\Box$

Now we can show a weak version of the mixed Bruhat decomposition for our groups $P_x$. (The weakness lies in the fact that we can not take two arbitrary points in $\Labar$ in the next statement.)
\begin{thm}
Let $x \in \Labar$ and $y \in \La$. Then we have $G = P_x N P_y$.
\end{thm}
{\bf Proof: } Let $\La_I$ be the component of $\Labar$ containing $x$. Then we can write $x = \sum_{i \in I} x_i \eta_i^I$ with some real coefficients $x_i$. We define a sequence of points in $\La$ by
\[ z_k = \sum_{i \in I} x_i \eta_i - k \sum_{ i \notin I} \eta_i.\]
Obviously, $z_k$ converges towards $x$. Now we choose a linear ordering $\prec$ on the set $\n$ in such a way that $i \in I$ and $j \notin I$ implies $i \prec j$. The set $\Phi^+ = \{ a_{ij} \in \Phi: i \prec j\}$ defines an order corresponding to some chamber. Note that for any root $a_{ij}$ in $\Phi^-$ we have
\[ a_{ij}(z_k) = \left\{ \begin{array}{ll} x_i -x_j&  \mbox{if } i,j \in I \\
-k-x_j & \mbox{if } i \notin I, j \in I \\
0 & \mbox{if } i,j \notin I.\end{array} \right. \]
Hence $a_{ij}(z_k)$ is bounded from above by a constant $c$ independent of $k$ and of the root $a_{ij} \in \Phi^-$. Therefore $U_{a_{ij}, z_k}$ is contained in $U_{a_{ij}, -c} = \{ u \in U_{a_{ij}}: \psi_a(u) \geq -c\}$, which is a compact subgroup of $U_{a_{ij}}$.
By 4.7, we find that all $U^-_{z_k}$ are contained in a compact subset of $U_{\Phi^-}$. 

 We have the ``usual'' mixed Bruhat decomposition for two points in $\La$ (see \cite{la}, 12.10), hence $G = P_{z_k} N P_y$ for all $k$. Using 4.8, we can  write an element $g \in G $ as
\[g = u_k^- u_k^+ n_k v_k\]
with $u_k^\pm \in U_{z_k}^\pm$, $n_k \in N$ and $v_k \in U_y$. Let us denote the kernel of the map $\nu: T \rightarrow \La$ by $Z \subset T$. Then the group $U_y^\wedge = U_y Z$ is compact and open in $G$ by \cite{la}, 12.12.
Since all $v_k$ lie in this compact subset, by switching to a subsequence we can assume that $v_k$ converges to some element $v \in U_y^\wedge$. Hence the sequence $u_k^- u_k^+ n_k$ is also convergent in $G$. Since the Weyl group $W = N/T$ is finite, by passing to a subsequence we can assume that $n_k = t_k n$ for some $t_k \in T$ and some fixed $n \in N$. Besides, we can assume that $u_k^-$ converges to some $u^- \in U_{\Phi^-}$, since the $u_k^-$ are contained in a compact subset of $U_{\Phi^-}$. Hence the sequence $u_k^+ t_k$ converges
in $G$. Its limit must be contained in the Borel group $U_{\Phi^+} T$. Hence $u_k^+$ converges towards some $u^+ \in U_{\Phi^+}$, and $t_k$ converges towards some $t \in T$. 

Using 4.7, $u_k^+$ is a product of $u_{a,k} \in U_{a,z_k}$ for all $a \in \Phi^+$, and, applying 4.3, we deduce that the $u_{a,k}$ converge towards some element $u_a \in U_{a,x}$. Hence we see that 
$u^+$ is contained in $U_{\Phi^+} \cap U_x = U_x^+ $. Similarly, $u^-$ lies in $U_x^-$. Therefore $g = u_k^- u_k^+ t_k n v_k$ converges towards $u^- u^+ t n v$, which is contained in $U_x^- U_x^+ N U_y \subseteq P_x N P_y$. Hence  $g$ lies indeed in $P_x N P_y$.\hfill$\Box$

Recall that $Z\subset T$ denotes the kernel of the map $\nu: T \rightarrow \La$, and that the group $U_0^\wedge = U_0 Z$ is compact. We define our compactification $\ov{X}$ of the building $X$ as 
\[\ov{X} = U_0^\wedge \times \Labar/\sim,\]
where the equivalence relation $\sim$ is defined as follows:
\begin{eqnarray*}
(g,x) \sim (h,y) & \mbox{iff  there exists an element } n \in N \\
~ & \mbox{such that } \nu(n)x= y \mbox{ and } g\inv h n \in P_x.\end{eqnarray*}
(Using 4.4, it is easy to check that $\sim$ is indeed an equivalence relation.)
We equip $\ov{X}$ with the quotient topology. The inclusion $U_0^\wedge \times \La \hookrightarrow G \times \La$ induces a bijection $(U_0^\wedge \times \La) / \sim \; \rightarrow X$, which is a homeomorphism if we endow the left hand side with the quotient topology (see \cite{bose}, p.221). Hence $X$ is open and dense in $\ov{X}$. 

We have a natural action of $U_0^\wedge$ on $\ov{X}$ via left multiplication on the first factor, which can be continued to an action of $G$ in the following way: If $g \in G$ and $(v,x) \in U_0^\wedge \times \Labar$, we can use the mixed Bruhat decomposition to write $gv = unh$ for some $u \in U_0^\wedge$, $n \in N$ and $h \in P_x$. Then we define $g(v,x) = (u, \nu(n) x)$. Using 4.4, one can show that this induces a well-defined action on $\ov{X}$. 

Mapping $x$ to the class of $(1,x)$, defines a map $\Labar \rightarrow \ov{X}$. This is injective, since by 4.8 we have $P_x  \cap N = N_x$. 
The $G$-action on $\ov{X}$ continues the $N$-action on $\Labar$, so that we will write $nx$ instead of $\nu(n) x$ for $x \in X$.

The following important fact follows immediately from the definition of $\ov{X}$: 
\begin{lem}
For all $x \in \Labar$ the group $P_x$ is the stabilizer of $x$ in $G$.
\end{lem}

We can use the mixed Bruhat decomposition to prove the following important fact:

\begin{prop}
For any two points $x \in \ov{X}$ and $y \in X$  there exists a compactified appartment containing $x$ and $y$,  i.e. there exists some $g \in G$ such that $x$ and $y$ both lie in $g \Labar$.
\end{prop}
{\bf Proof: } We can assume that $y$ lies in $\La$. The point $x$ lies in $h \Labar$ for some $h \in G$, so $x = h x'$ for some $x' \in \Labar$. By our mixed Bruhat decomposition  4.9 we can write  $h= qnp$ for $q \in P_y$, $n \in N$ and $p \in P_{x'}$. Therefore $x = h x' = qnp x' = qn x' \in q \Labar$, and $y = qy$ lies also in $q \Labar$, whence our claim.\hfill$\Box$

\section{Properties of $\ov{X}$}
In this section we want to check that $\ov{X}$ is compact and we want to identify it with the set 
$\bigcup_{W \subseteq V} X(PGL(W))$. 

Hence we see that we can compactify the Bruhat-Tits building for $PGL(V)$ by attaching all the Bruhat-Tits building for $PGL$ of the smaller subspaces at infinity.

The following lemma is similar to \cite{la}, 8.11.

\begin{lem}
Let $z$ be a point in $\La_I$ and $y$ be a point in $\La_J$ for some $I \subseteq J\subseteq \n$. Then we find a chamber such that the corresponding set of positive roots $\Phi^+$ satisfies $U_y^+ \subseteq U_z^+$. 
\end{lem}
Note that the assumptions are fulfilled if $J = \n$, i.e. if $y$ lies in $\La$.

{\bf Proof: }
Since $I \subseteq J$, we can define a projection map $ \La_J \rightarrow \La_I$, which we also denote by $r_I$. To be precise, $r_I$ maps a point $\sum_{i \in J} x_i \eta_i^J \in \La_J$ to $\sum_{i \in I} x_i \eta_i^I \in \La_I$. Put $y^\ast = r_I(y)$.

Recall that we denote by $\Phi_I$ the set of roots of $\bt_I$ in $\bg_I = PGL(V_I)$. There exists a chamber in $\La_I^\ast$ with respect to $\Phi_I$ such that the corresponding subset $\Phi_I^+$ of positive roots satisfies
$U_{y^\ast}^{I+} \subseteq U_z^{I+}$.
Here $U_z^{I+}$ is defined exactly as the groups $U_x^+$ for $x \in \La$ in section 4, just replacing $\La$ by $\La_I$. 

Now $\Phi_I = \{ b_{ij}: i,j \in I\}$, where $b_{ij}$ is the character mapping a diagonal matrix with entries $t_l$ (for $l \in I$) to $t_i / t_j$. We define a linear ordering on $I$ by $i \prec j $ iff $b_{ij} \in \Phi^+_I$. This can be continued to a linear ordering on $\n$ in such a way that $i \prec j$ whenever $i \in I$ and $j \notin I$. Let us put $\Phi^+ = \{a_{ij}: i \prec j\}$. We claim that this satisfies our claim. In fact, take $a = a_{ij} \in \Phi^+$. If $i$ and $j$ are contained in $I$, we use 4.2 and the construction of $\Phi^+_I$ to deduce that $U_{a,y} \subseteq U_{a,z}$. By definition of $\Phi^+$ it cannot happen that $j$ is in $I$, but $i$ is not. Hence the only remaining case is that $j \notin I$. Then $f_z(a) = - \infty$ by 4.1, so that trivially $U_{a,y} \subseteq U_{a,z}$.\hfill$\Box$

\begin{thm}
Fix a nonempty $\Omega \subset \Labar$. We denote by $(\ast)$ the following condition:
\begin{eqnarray*} (\ast) \quad  &\mbox{}  \mbox{The set } \{J \subseteq \n : \Omega \cap \La_J \neq \emptyset  \}
\mbox{ contains a maximal element}\\ ~ & \mbox{} \mbox{with respect to inclusion.}
\end{eqnarray*}
If $(\ast)$ if satisfied, then $P_\om = \cap_{x \in \om} P_x$. In particular, $P_\om$ is the stabilizer of $\om$.
\end{thm}
Note that $(\ast)$ is satisfied if $\om \cap \La$ is not empty. 

{\bf Proof: }To begin with, the inclusion $P_\om \subseteq P_x$ for all $x \in \om$ is trivial. 
Let $J$ be the maximal subset of $\n$ satisfying $\om \cap \La_J \neq \emptyset$, and choose some $x_0 \in \om \cap \La_J$. We will first prove that $P_{\om^\sim} \cap P_x = P_{\om^\sim \cup \{x\}}$ for all $x \in \om$ and all subsets $\om^\sim$ of $\om$ containing $x_0$. Assume that $x$ lies in in $\La_I$ for some $I \subseteq J$. By 5.1, we find some set of positive roots $\Phi^+$ such that $U_{x_0}^+ \subseteq U_x^+$, hence also $U_{\om^\sim}^+ \subseteq U_{x}^+$. Let now $g$ be an element in $P_{\om^\sim}\cap P_x$, and write $g = n u^- u^+$ for some $n \in N_{\om^\sim}$ and some $u^\pm \in U_{\om^\sim}^\pm$ (using 4.8). Since  $g$ and $u^+$ are contained in $P_x$, this also holds for $n u^-$, so that we can write $n u^- = m v^+ v^-$ for some $m \in N_x$ and $v^\pm \in U_x^\pm$. So $m\inv n$ is contained in $N \cap U_{\Phi^+} U_{\Phi^-}$, hence trivial by \cite{boti}, 5.15. We find that $m = n$ and $u^- = v^+ v^-$, which implies $v^+ = 1$. Hence $n$ is contained in $N_x \cap N_{\om^\sim} = N_{\om^\sim \cup \{x\}}$, and $u^-$ is contained in $U_x^- \cap U_{\om^\sim}^-$. 
Note that for all $a \in \Phi$ we have the inclusion $U_{a,\om^\sim} \cap U_{a,x} \subseteq U_{a,\om^\sim \cup \{x\}}$, so that we can use 4.7 to deduce $u^- \in U_{\om^\sim \cup \{x\}}^-$. 
A similar argument as above gives $u^+ \in  U_{\om^\sim}^+ \cap U_x^+ \subset U^+_{\om^\sim \cup \{x\}}$.
Hence our claim is proven.

Therefore any finite subset $\om^\sim \subset \om$ containing $x_0$ satisfies
our claim, i.e.  $P_{\om^\sim} = \cap_{x \in \om^\sim} P_x$.

We can write $\om = \bigcup_{\sigma \in \Sigma} \om_\sigma$, where $\om_\sigma$ for $\sigma \in \Sigma$  runs over all finite subsets of $\om$ containing $x_0$. Let us consider some  $g \in \bigcap_{x \in \om} P_x = \bigcap_{\sigma \in \Sigma} P_{\om_\sigma}$. We fix some set of positive roots $\Phi^+$ and write $g = n_\sigma u_\sigma^+ u_\sigma^-$ for $n_\sigma \in N_{\om_\sigma}$ and $u_\sigma^\pm \in U_{\om_\sigma}^\pm$ by 4.8. Put $T_{x_0} = T \cap N_{x_0}$. Then $N_{x_0} / T_{x_0}$ is finite. By the pidgeon hole principle, there must be one class $m$ in $N_{x_0} / T_{x_0}$ such that the set $\Sigma'$ of all the $\sigma \in \Sigma$ so that  $n_{\sigma}$ is equal to $m$ modulo $T_{x_0}$ still has the property that $\bigcup_{\sigma \in \Sigma'} \om_{\sigma} = \om$. (If not, we could find for any class in $N_{x_0} / T_{x_0}$ an element in $\om$ not contained in any $\om_{\sigma}$ with the property that $n_\sigma$ lies in our class. Collecting these elements together with $x_0$ in some finite set gives a contradiction.)

Hence for $\sigma \in \Sigma'$ we can write $n_\sigma = m t_\sigma$ for some fixed $m \in N_{x_0}$ and some $t_\sigma \in T_{x_0}$. For $\sigma$ and $\tau$ in $\Sigma'$ we get $t_\sigma u_\sigma^+ u_\sigma^- = t_\tau u_\tau^+ u_\tau^-$. Using the fact that $T$ normalizes $U_{\Phi^+}$ and $U_{\Phi^-}$ and that  $N \cap U_{\Phi^+} U_{\Phi^-}$ is trivial, we find $u_\sigma^- = u_\tau^-$, $u_\sigma^+ = u_\tau^+$ and $t_\sigma = t_\tau$. Therefore the elements $t = t_\sigma$, $u^\pm = u_\sigma^\pm$ are independent of the choice of $\sigma \in \Sigma'$. Note that by definition of the $T$-action on $\Labar$ the element $t$ in $T_{x_0}$ stabilizes not only $x_0$, but also every point in the components $\La_I$ for $I \subseteq J$. Thus $t \in N_\om$. Besides, by 4.7 we deduce $u^+ \in \bigcap_{\sigma \in \Sigma'} U_{\om_\sigma}^+ = U_\om^+$. Similarly, $u^- \in U_\om^-$. Hence we find indeed that $g = mtu^+ u^-$ is contained in $P_\om$.\hfill$\Box$

\begin{cor}
Let $\om$ and  $\om'$ be two non-empty subsets of $\Labar$ such that $\om$, $\om'$ and $\om \cup \om'$ satisfy condition $(\ast)$ in the Theorem. Then
$P_\om \cap P_{\om'} = P_{\om \cup \om'}$.
\end{cor}
{\bf Proof: } This is an immediate consequence from 5.2.\hfill$\Box$

The following result is similar to \cite{la}, 9.6.
\begin{prop}
Let $g\in G$ and let $J$ be a subset of $\n$ so that $\Lambda_J \cap g\inv \Labar$ is not empty.  Then there exists some element $n \in N$ such that
\[ g x = n x \mbox{ for all } x \in g\inv \Labar \cap (\bigcup_{I \subseteq J} \La_I).\]
\end{prop}
{\bf Proof: }Note that the set $\om =  g\inv \Labar \cap (\bigcup_{I \subseteq J} \La_I)$ satisfies condition $(\ast)$ from 5.2. Fix some $x_0 \in \La_J \cap g\inv \Labar$. For all $x \in \om$ we have $g\inv N \cap P_x \neq \emptyset$, since $x = g\inv y$ for some $y \in \Labar$. 

We will now show  that for all finite subsets $\Delta$ of $\om$ containing $x_0$ we have $g\inv N \cap P_{\Delta} \neq \emptyset$. Let us suppose that this claim holds for some $\Delta$ and let us show it for $\Delta \cup \{x\}$, where $x$ is some point in $\om$. So there is some $n_\Delta \in N$ with $g\inv n_\Delta \in P_{\Delta}$. We also find some $n_x \in N$ satisfying $g\inv n_x \in P_x$. 
By 5.1, we find a set of positive roots $\Phi^+$ such that $U_{x_0}^+ \subseteq U_x^+$, so that also $U_\Delta^+ \subseteq U_x^+$. Hence we apply 4.8 to deduce
\[n_\Delta\inv n_x \in P_{\Delta} P_x = N_\Delta U_\Delta^- U_\Delta^+ U_x^+ U_x^- N_x = N_\Delta U_\Delta^- U_x^- U_x^+ N_x \subseteq N_\Delta U_{\Phi^-} U_{\Phi^+} N_x.\]
Since $N \cap U_{\Phi^-} U_{\Phi^+}$ is trivial, we find  $n_\Delta' \in N_\Delta$ and $n_x' \in N_x$ such that $n = n_x n_x' = n_\Delta n_\Delta'$ satisfies $g\inv n \in P_\Delta \cap P_x$,  which is equal to $P_{\Delta \cup \{x\}}$ by 5.3. This proves our claim. 

Now we write as in Theorem 5.2 $\om = \bigcup_{\sigma \in \Sigma} \om_\sigma$, where $\om_\sigma$ runs over all finite subsets of $\om$ containing $x_0$. For all $\sigma$ we choose some $n_\sigma \in N$ such that $g\inv n_\sigma \in P_{\om_\sigma}$. Put $n_\sigma = n_0$ if $\om_\sigma$ is the set $\{x_0\}$. The same argument as in 5.2 shows that we can find a subset $\Sigma' \subseteq \Sigma$ such that $\om = \bigcup_{\sigma \in \Sigma'}\om_\sigma$ and such that $n_\sigma\inv n_0$ is equal to some fixed $m \in N_{x_0}$ modulo $T_{x_0}$ for all $\sigma \in \Sigma'$. Since any element in $T_{x_0}$ leaves the components $\La_I$ for $I \subseteq J$ pointwise invariant, it also stabilizes $\om$. Therefore $n_\sigma\inv n_0 m\inv$ lies in $N_\om \subset P_\om$. Since $g\inv n_\sigma$ is contained in $P_{\om_\sigma}$, the same holds for $g\inv n_0 m\inv$, so that $g\inv n_0 m\inv$ lies in $\bigcap_{\sigma \in \Sigma'}  P_{\om_\sigma}$, which is equal to $P_\om$ by 5.3. Hence $n = n_0 m\inv$ satisfies $nx = gx$ for all $x \in \om$, as desired.\hfill$\Box$

Now we can prove
\begin{thm}
$\ov{X}$ is compact.
\end{thm}
By \cite{bout}, I.10.4, Proposition 8, we know that since $U_0^\wedge \times \Labar$ is compact, the quotient after $\sim$ is Hausdorff (hence compact), iff the relation $\sim$ is closed in $(U_0^\wedge \times \Labar) \times (U_0^\wedge \times \Labar)$. 

Let $(u_k)_k$ and $(v_k)_k$ be sequences in $U_0^\wedge$ converging to $u$ respectively $v$, and let $x_k$ and $y_k$ be sequences in $\Labar$ converging to $x$ respectively $y$, so that $(u_k, x_k) \sim (v_k, y_k)$. We have to show that $(u,x) \sim (v,y)$. By definition of $\sim$ we have $x_k = u_k\inv v_k y_k$, so that both $0$ and $x_k$ lie in $\Labar$ and in $u_k\inv v_k \Labar$. Using 5.4, we find some $n_k \in N$ such that $n_k z = v_k\inv u_k z$ for all $z \in \Labar \cap u_k\inv v_k \Labar$. In particular, $n_k$ lies in $N_0$ and $n_k x_k = v_k\inv u_k x_k = y_k$. 

Hence $g_k = u_k\inv v_k n_k$ lies in $P_{x_k} \cap P_0 = P_{\{0, x_k\}}$. Since $N_0$ is compact, we can pass to a subsequence and assume that $n_k$ converges towards some $n \in N_0$, so that $g_k$ converges towards some $u\inv v n$. 

By 4.8, we can write $g_k = w_k^- w_k^+ m_k$ for some $w_k^\pm \in U_{\{0, x_k\}}^\pm$ and $m_k \in N_{\{0, x_k\}}$. We can again assume that $m_k$ converges towards some $m \in N_0$. Since $N$ acts continuously on $\Labar$, $m$ lies also in $N_x$. Besides, $U_0^-$ is compact, so that we can assume that $w_k^-$ converges towards some $w^- \in U_0^-$. Using 4.7 and 4.3 we find that $w^-$ lies in fact in $U_x^-$. Now $w_k^+$ also converges towards some $w^+$ which lies in $U_x^+$ by the same argument. Therefore $u\inv v n$ lies in $P_x$. Since $N$ acts continuously on $\Labar$, we have $nx = y$, so that indeed $(u,x) \sim (v,y)$.\hfill$\Box$

\begin{thm}
The space $\ov{X}$ is contractible.
\end{thm}
{\bf Proof: }Recall from 3.4 that $\Labar$ is contractible. If $x = \sum_{ i \neq j \in J} x_j (-\eta_j^J)$ is a point in $E_i \cap \La_J$, where $J$ contains $i$, then the contraction map is given by
\[r(x,t) = \left\{ \begin{array}{ll} x, & \mbox{ if } t = 0 \\
\sum_{i \neq j \in J }\frac{(1-t) x_j}{1+ t x_j} (-\eta_j) + \sum_{ j \notin I} \frac{1-t}{t} (- \eta_j), & \mbox{ if }t \neq 0. \end{array} \right.\]
Now we define
\begin{eqnarray*} R: & U_0^\wedge \times \Labar \times [0,1] & \longrightarrow U_0^\wedge \times \Labar \\
~ & ((g,x),t) & \longmapsto (g, r(x,t)).
\end{eqnarray*}
Obviously, $R$ is continuous. In order to show that $R$ is a contraction map for $\ov{X}$, it suffices to prove that it is compatible with our equivalence relation. 

Let us first check that $n r(x,t) = r(nx, t)$ for all $n \in N_0$ and $x \in \Labar$. Write $n = tp$ for $t \in T$ and a permutation matrix $p$. Then $t$ lies in $N_0$, so that it acts trivially on all points in $\Labar$. Besides, a straightforward calculation shows $p r(x,t) = r (px,t)$. 

Now assume that $(g,x)$ and $(h,y)$ in $U_0^\wedge \times \Labar$ are equivalent.
Hence there is some $n \in N$ with $nx =y$ and $g\inv h n \in P_x$. Using 5.4, we can assume that $n$ lies in fact in $N_0$.

Now fix some $t \in [0,1]$. We already know that $r(y,t) = r(nx,t) = n r(x,t)$, so that 
our claim, namely $(g, r(x,t)) \sim ( h, r(y,t))$ is proven, if we show that $g\inv h n \in P_{r(x,t)}$. We already know that $g\inv h n$ lies in $P_0 \cap P_x$, which is equal to $P_{\{0,x\}}$ by 5.3.
Let us put $x_t = r(x,t)$. 
We are done if we show that $P_{\{0,x\}}\subset P_{x_t}$. We can assume that $t >0$. The point $x$ is contained in some $E_i \cap \La_J$. Recall from 4.8 that $P_{\{0,x\}} = U_{\{0,x\}}^- U_{\{0,x\}}^+ N_{\{0,x\}}$ and $P_{x_t} = U_{x_t}^- U_{x_t}^+ N_{x_t}$ for some fixed $\Phi^+$.
A straightforward calculation yields
\begin{eqnarray*}
f_{x_t}(a_{kl})  \leq  & f_x(a_{kl}),& \mbox{ if } k\neq i \mbox{ and } l \mbox{ are in } J \mbox{ and }f_x(a_{kl}) \geq 0, \mbox{ or if }k \notin J \mbox{ and } l \in J\\
f_{x_t}(a_{kl})  \leq & 0,& \mbox{if } k\neq i \mbox{ and } l \mbox{ are in } J \mbox{ and } f_x(a_{kl}) < 0, \\
~ & ~ & \mbox{or if } k = i \mbox { and } l \in J, \mbox{ or if } l \notin J
\end{eqnarray*}
Hence for all $a \in \Phi$ we have $U_{a,0} \subset U_{a, x_t}$ or $U_{a,x}\subset U_{a, x_t}$, which implies that $U_{\{0,x\}}^+ \subset U_{x_t}^+$ and $U_{\{0,x\}}^- \subset U_{x_t}^-$. Besides, we have $N_{\{0,x\}} \subset N_{x_t}$, since $r$ is compatible with the action of $N_0$, so that our claim follows.\hfill$\Box$

The following result shows that the boundary of our compactification $\ov{X}$ consists of the Bruhat-Tits buildings of all groups $PGL(W)$, where $W$ is a non-trivial subspace of $V$:

\begin{thm}
There is a bijection between $\ov{X}$ and the set \[\bigcup_{0 \neq W \subseteq V} PGL(W) = X \cup \bigcup_{0 \neq W \subset V} PGL(W).\]
\end{thm}

{\bf Proof: }Let us first fix some non-empty $I \subset \n$. We will start by embedding the building corresponding to $\bg_I = PGL(V_I)$ in $\ov{X}$. Recall that we write $\bg^{V_I}$ for the subgroup of $\bg$ leaving $V_I$ invariant, and $\rho_I$ for the map $\bg^{V_I} \rightarrow \bg_I$. Then the building for $\bg_I$ is defined as $G_I \times \La_I / \sim$, where $\sim$ is the equivalence relation from section 2 (replacing $V$ by $V_I$ everywhere). Let $(h,x)$ be an element of $G_I \times \La_I$. We choose an arbitrary lift $h^\uparrow$ of $h$ in $G^{V_I}$, and map $(h,x)$ to the point $h^\uparrow x$ in $\ov{X}$. This induces a map 
\[ j_I: G_I \times \La_I \longrightarrow \ov{X}.\]
We claim that this is independent of the choice of a lift. We have to show that any element $g$ in the kernel of $\rho_I$ stabilizes each $x \in \La_I$. Note that $G^{V_I}$ is a parabolic subgroup of $G$. Let us fix a Borel group $B \supset T$ contained in $G^{V_I}$. Then there is a set of positive roots $\Phi^+$ in $\Phi$ such that $B = U_{\Phi^+} T$. Besides, we can write $G^{V_I} = B W' B$ for some subgroup $W' $ of $W$. Hence we find that $G^{V_I} = U_{\Phi^+} N^{V_I} U_{\Phi^+}$, where $N^{V_I} = N \cap G^{V_I}$. Note that since $U_{\Phi^+}$ is contained in $G^{V_I}$, no root of the form $a = a_{ij}$ such that $j \in I$, but $i \notin I$ can be contained in $\Phi^+$. We write $g = u^+ n v^+$ for some $n $ in $N^{V_I}$ and $u^+, v^+$ in $U_{\Phi^+}$. Using $\rho_I(g) = 1$, we find that $\rho_I(n) = 1$ and $\rho_I(u^+ v^+) = 1$. Hence $n$ stabilizes each $x \in \La_I$, so $n \in P_x \cap N = N_x$.

Let us write $[I]$ for the subset of all roots that are linear combinations of roots $a_{ij}$ with $i$ and $j$ in $I$. We write $u^+ = u_1^+ u_2^+$ and $v^+ = v_2^+ v_1^+$ for some $u_1^+, v_1^+ \in m(\prod_{a \in \Phi^+ \backslash [I]} U_a)$ and $u_2^+, v_2^+ \in m(\prod_{a \in \Phi^+ \cap [I]} U_a)$, where $m$ is the multiplication map. As in 4.2, we find that $\rho_I(u_1^+)$ and $\rho_I(v_1^+)$ are trivial, and that $\rho_I$ is injective on all $U_a$ for $a \in \Phi^+ \cap [I]$. Hence we deduce that $u_2^+ v_2^+ $ is trivial. Besides, $u_2^+$ commutes with $n$, so that $g = u_1^+ u_2^+ n v_2^+ v_1^+ = u_1^+ n v_1^+$. Recall that the roots $a_{ij}$ for $i \notin I$ and $j \in I$ do not lie in $\Phi^+$, so that by 4.1 we have $U_a = U_{a,x}$ for all $a \in \Phi^+\backslash [I]$ and all $x \in \La_I$. Hence $g$ is indeed contained in $P_x$ for all these $x$. 

Now we claim that $j_I$ induces an injection
\[j_I:  X(G_I) = G_I \times \La_I / \sim \; \;\hookrightarrow \ov{X}.\]
It suffices to show that for all $x \in \La_I$ the map $\rho_I$ induces a surjection 
\[\rho_I: P_x \longrightarrow P_x^I,\]
where $P_x^I$ is defined in the same way as $P_x$ (see section 2), just replacing $V$ by $V_I$. Note that this also implies $P_x = \rho_I\inv(P_x^I)$. 

First of all, let us show that $P_x$ is  contained in $G^{V_I}$, so that our statement makes sense. By 4.2, we find that $U_x$ is contained in $G^{V_I}$. Besides, if $n$ is an element of $N_x$, we find that $n= t p$ for a torus element $t$ and a permutation matrix $p$. The corresponding permutation must leave $I$ intact, so that indeed $n \in N^{V_I}$. 

Besides, by 4.2 we know that $\rho_I$ maps $U_x$ to $U_x^I$, and it is easy to see that $\rho_I$ also maps $N_x$ to $N_x^I$, so that  we get a homomorphism $ \rho_I: P_x \rightarrow P_x^I$. Now take some $h = n u $ in $P_x^I = N_x^I U_x^I$, and choose a lift $n^\uparrow$ of $n$ in $N^{V_I}$. Since $n$ stabilizes $x$ (viewed in the appartment $\La_I$), the lift $n^\uparrow$ stabilizes the point $x \in \Labar$. Hence $n^\uparrow$ lies in $P_x$. By 4.2, the element $u$ has a lift $u^\uparrow$ in $U_x$, hence $h^\uparrow = n^\uparrow u^\uparrow$ is an element in $P_x$ projecting to $h$ via $\rho_I$. This proves surjectivity. 

Let now $W$ be an arbitrary non-trivial subspace of $V$. Then there is a linear isomorphism 
\[f: V_I \longrightarrow W\]
for some $I \subset \n$. Let $\bf{S}$ be the maximal torus in $PGL(W)$ induced by the diagonal matrices with respect to the basis $f(v_i)$ (for all $i \in I$). Conjugation by $f$ induces an isomorphism $PGL(W) \rightarrow PGL(V_I) = \bg_I $, which maps ${\bf S}$ to ${\bf T}_I$, and the normalizer $N({\bf S})$ of ${\bf S}$ to $N({\bf T}_I)$, the normalizer of ${\bf T}_I$ in $G_I$. Hence we get an $\R$-linear isomorphism
\[\tau: X_\ast({\bf S})_\R \iso X_\ast({\bf T}_I)_\R = \La_I.\]
One can check that for all $x \in X_\ast({\bf S})_\R$ and $n \in N({\bf S})$ we have $\tau(nx) = (f\inv n f) \tau(x)$. 

Choose some $f^\uparrow \in G$ whose restriction to $V_I$ is given by $f$. Then we  define a map
\[j_W: PGL(W) \times X_\ast({\bf S})_\R \longrightarrow G_I \times \La_I \stackrel{j_I}{\longrightarrow} \ov{X} \stackrel{f^\uparrow}{\longrightarrow} \ov{X},\]
where the first map is given by
\[ (g,x) \longmapsto (f\inv g f, \tau(x)).\]
Since this maps the equivalence relation on $PGL(W) \times X_\ast({\bf S})_\R$ defining the building $X(PGL(W))$ to the equivalence relation on $G_I \times \La_I$ defining the building $X(G_I)$, we have an injection
\[j_W: X(PGL(W)) \iso X(G_I) \stackrel{j_I}{\longrightarrow} \ov{X} \stackrel{f^\uparrow}{\longrightarrow} \ov{X}.\]
Of course, we have to check that this is well-defined. First of all, if $f$ is fixed, then $j_W$ does not depend on the choice of a lifting $f^\uparrow$, since two such liftings differ by something in the kernel of $\rho_I$, and this acts trivially on the image of $j_I$, as we have seen above.

What happens if we choose another isomorphism $g: V_J \rightarrow W$ for some $J \subset \n$? First let us consider the case that we construct $j_W$ using $f' = s \circ f$ for some isomorphism $s : W \rightarrow W$. Then we also use a different construction of the building $X(PGL(W))$, since we use another torus. Following \cite{la}, 13.18, we find that there exists a unique $PGL(W)$--equivariant isometry between these two constructions (which we tacitly use to identify them). A straightforward calculation shows that our map $j_W$ is compatible with this isometry. 

Now assume that we use an isomorphism $g : V_J \rightarrow W$ to construct $j_W$. Then there exists some isomorphism $r : V_I \rightarrow V_J$ mapping the basis $v_i$ for $i \in I$ to the basis $v_j$ for $j \in J$ in some way. The map $r$ can be lifted to a permutation matrix $n \in N$. It is easy to see that this implies that our construction of $j_W$ does indeed not change if we choose $ g\circ r$ instead of $g$. Hence $j_W$ is well-defined. 

We put all these maps $j_W$ together and get a map
\[ j : X \cup \bigcup_{0 \neq W \subset V} X(PGL(W)) \longrightarrow \ov{X},\]
which is obviously surjective. It remains to show that $j$ is injective.
Let us first assume that we have some $(g,x) \in PGL(W) \times X_\ast({\bf S})_\R$ and $(h,y) \in PGL(W') \times X_\ast({\bf S}')_\R$ such that $j_W(g,x) = j_{W'}(h,y)$. Hence there are isomorphisms $f_I: V_I \rightarrow W$ and $f_J: V_J \rightarrow W'$ and points $x_0 \in \La_I$, $y_0 \in \La_J$ such that $g^\uparrow f_I^\uparrow x_0 = h^\uparrow f_J^\uparrow y_0$. In particular, there is some $n \in N$ mapping $x_0$ to $y_0$. Hence $n$ maps $\La_I$ to $\La_J$, which implies $n V_I = V_J$. We have already seen that $P_{x_0} \subset G^{V_I}$, so that    $h^\uparrow f_J^\uparrow n V_I = g^\uparrow f_I^\uparrow V_I$, which implies $W = f_I(V_I) = f_J (V_J) = W'$. Since we already know that $j_W $ is injective, our claim follows.\hfill$\Box$

To conclude this paper, let us show that we can identify the vertices in $\ov{X}$ with the equivalence classes $\{N\}$ of $R$-modules of arbitrary rank in $V$. Together with Proposition 3.1, this is the link to Mustafin's paper \cite{mu}.

We call a point $x$ in $\Labar$, say $x \in \La_I$,  a vertex in $\Labar$, if it is a vertex in the appartment $\La_I$. By section 2 this means, that $x = \sum_{i \in I} x_i \eta_i^I$ with integer coefficients $x_i$. We call a point $y$ in $\ov{X}$ a vertex if $y = gx$ for some $g \in G$ and some vertex $x \in \Labar$, and we denote the set of vertices in $\ov{X}$ by $\ov{X}^0$. 

We call two $R$-lattices in $V$ equivalent, if they differ by a factor in $K\tim$.
Let us denote by $\ov{\cL}$ the set of all equivalence classes of $R$-lattices in $V$ of {\em{arbitrary}} positive rank. We write $\{M\}$ for the class of such a lattice.

Our last result shows that the vertices in $\ov{X}$ correspond to elements of $\ov{\cL}$, i.e. lattice classes of arbitrary rank, which explains the title of this paper. 
\begin{lem}
The $G$-equivariant bijection $\varphi: \cL \longrightarrow X^0$ can be continued to a $G$-equivariant bijection 
\[ \varphi: \ov{\cL} \longrightarrow \ov{X}^0\]
in the following way: We write $\{M\} \in \ov{\cL}$ as $\{M\} = g \{L\}$, where 
$L = \sum_{i \in I} \pi^{k_i} R v_i$ for some non-empty $I \subseteq \n$. Then 
\[ \varphi(\{M\}) = g ( \sum_{i \in I} k_i (-\eta_i^I)).\]
\end{lem} 
{\bf Proof: }We only need to check that $\varphi$ is injective and well-defined, which amounts to the following claim: For the vertex $x = \sum_{i \in I} k_i (-\eta_i^I) \in \La_I$ let $L_x =   \sum_{i \in I} \pi^{k_i} R v_i$. Then $P_x$ is the stabilizer of $\{L_x\}$. Using 4.1, it is easy to see that all $U_{a,x}$ leave $\{L_x\}$ invariant. Besides, $N_x$ leaves $\{L_x\}$ invariant, so that $P_x$ is contained in $S_{\{L_x\}}$, the stabilizer of $\{L_x\}$. 

Let now $g$ be an element in $S_{\{L_x\}}$. Note that $L_x$ is a lattice of full rank in $V_I$, so that $g$ is contained in $G^{V_I}$, and $\rho_I(g) \in PGL(V_I)= G_I$ stabilizes $\{L_x\}$. In $G_I$ we have the (usual) Bruhat decomposition 
$\rho_I(g) = pnq$ with $p$ and $q$ in $P_x^I$ and $n \in N(T_I)$ (see \cite{la}, 12.10). Here $P_x^I$ is the stabilizer of $x$ in the building for $G_I$, so that $p$ and $q$ leave the class $\{L_x\}$ in $V_I$ invariant. Hence $n \in G_I$ also stabilizes $\{L_x\}$. A straightforward calculation now shows that $n$ fixes $x \in \La_I$. Therefore $\rho_I(g)$ lies in $P_x^I$, so that $g$ is indeed contained in $P_x$, as we have seen in the proof of 5.7.\hfill$\Box$

Of course, one could also set up a bijection between  $\ov{\cL}$ and $\bigcup_{0 \neq W \subseteq V} X(PGL(W))^0$ by using analogues of the map $\varphi$ on every single building $X(PGL(W)$. It is easy to see that this is compatible with the map we just described, und the identification of $\ov{X}$ and $\bigcup_{0 \neq W \subseteq V} X(PGL(W))$ given in 5.7.

\end{document}